\documentclass[12pt,a4paper]{article}
\usepackage{mathrsfs}
\usepackage{color}
\usepackage{amsbsy}
\usepackage{amsmath}
\usepackage{amssymb}
\usepackage{multicol}
\usepackage{color,verbatim,lscape}
\usepackage{graphicx}
\usepackage{caption,subcaption}
\usepackage[round]{natbib}
\usepackage{pdflscape}
\usepackage{hyperref}
\usepackage{float}
\usepackage{array}
\usepackage{makecell}
\usepackage{soul}
\usepackage{enumitem}
\hypersetup{
	colorlinks=true,
	citecolor=blue,
	linkcolor=blue,
	urlcolor=blue,
	breaklinks=true
}
\usepackage[nolabel, final]{showlabels}
\usepackage{booktabs}
\usepackage{multirow}
\usepackage{float}
\usepackage{rotfloat}

\allowdisplaybreaks \textheight 9.0 in \textwidth 6.6 in
\makeatletter 

\usepackage[linesnumbered, ruled]{algorithm2e}

\def\Var{\mbox{Var}}
\def\Cov{\mbox{Cov}}
\def\Exp{\mbox{E}}
\def\SSR{\mbox{SSR}}

\def\CV{\mbox{CV}}

\def\reg{\mathrm{reg}}
\def\pf{\mathrm{pf}}
\def\eig{\mathrm{eig}}
\def\ma{\mathrm{ma}}
\def\cv{\mathrm{cv}}

\newcommand{\F}{{\bf F}}

\newcommand{\A}{{\bf A}}
\newcommand{\w}{{\bf w}}

\newcommand{\bbW}{\mathbb{W}}
\newcommand{\bbA}{\mathbb{A}}
\newcommand{\bbB}{\mathbb{B}}
\newcommand{\M}{{\bf M}}

\newcommand{\s}{{ (s)}}
\newcommand{\x}{{\bf x}}

\newcommand{\bP}{{\bf P}}
\newcommand{\J}{\mathbf{J}}
\newcommand{\I}{{\bf I}}
\newcommand{\y}{{\bf y}}

\newcommand{\Q}{{\bf Q}}
\newcommand{\calR}{{\mathcal{R}}}
\newcommand{\calA}{{\mathcal{A}}}
\newcommand{\calB}{{\mathcal{B}}}
\newcommand{\calC}{{\mathcal{C}}}
\newcommand{\calD}{{\mathcal{D}}}
\newcommand{\calE}{{\mathcal{E}}}

\newcommand{\calX}{{\mathcal{X}}}
\newcommand{\calY}{{\mathcal{Y}}}
\newcommand{\calI}{{\mathcal{I}}}

\newcommand{\f}{{\bf f}}
\renewcommand{\b}{\mathbf{b}}

\newcommand{\e}{{\bf e}}
\newcommand{\bk}{{\bf k}}
\newcommand{\sums}{\sum_{s=1}^{S}}
\newcommand{\sumjs}{\sum_{j=1}^{S}}

\newcommand{\tp}{^\top}
\newcommand{\one}{\mathbf{1}}

\newcommand{\SAIC}{\textrm{SAIC}}
\newcommand{\SBIC}{\textrm{SBIC}}

\newcommand{\balpha}{\boldsymbol{\alpha}}
\newcommand{\bphi}{\boldsymbol{\phi}}
\newcommand{\bpsi}{\boldsymbol{\psi}}
\newcommand{\bSigma}{\boldsymbol{\Sigma}}

\newcommand{\bbeta}{\boldsymbol{\beta}}

\newcommand{\bmu}{\boldsymbol{\mu}}
\newcommand{\bnu}{\boldsymbol{\nu}}
\newcommand{\brho}{\boldsymbol{\rho}}

\newcommand{\Nor}{\mathcal{N}}

\newcommand{\bias}{\textrm{Bias}}

\newcommand{\model}{\textrm{model}}
\newcommand{\XIC}{\textrm{XIC}}

\DeclareMathOperator{\tr}{tr}
\DeclareMathOperator{\argmin}{argmin}
\DeclareMathOperator{\argmax}{argmax}
\DeclareMathOperator{\diag}{diag}

\newtheorem{proposition}{Proposition}

\newtheorem{remark}{Remark}

\newenvironment{proof}[1]{\par\textbf{\emph{Proof of #1}.}}{\par\hfill \ensuremath{\Box}\par}

\graphicspath{{figure/}}

\date{}
\title{A theoretical comparison of weight constraints in forecast combination and model averaging}

\author{Jiahui Zou\thanks{Capital University of Economics and Business. 
}\quad\qquad  Andrey Vasnev\thanks{University of Sydney. 
}\quad\qquad Wendun Wang\thanks{
Erasmus University Rotterdam; Tinbergen Institute.
}\quad\quad Xinyu Zhang\thanks{Chinese Academy of Sciences. 
}}
\date{\today}

\begin{document}

\maketitle


\textbf{Abstract:}
	Forecast combination and model averaging have become popular tools in forecasting and prediction, both of which combine a set of candidate estimates with certain weights and are often shown to outperform single estimates. A data-driven method to determine combination/averaging weights typically optimizes a criterion under certain weight constraints. While a large number of studies have been devoted to developing and comparing various weight choice criteria, the role of weight constraints on the properties of combination forecasts is relatively less understood, and the use of various constraints in practice is also rather arbitrary. In this study, we summarize prevalent weight constraints used in the literature, and theoretically and numerically compare how they influence the properties of the combined forecast. Our findings not only provide a comprehensive understanding on the role of various weight constraints but also practical guidance for empirical researchers how to choose relevant constraints based on prior information and targets. \\


\clearpage

\section{Introduction}\label{sec:intr}
Forecasting and prediction are among the most important tasks in economic analysis, in which forecast combination and model averaging techniques have gained increasingly popularity and even become benchmark methods in some contexts since the seminal work by \citet{Bates1969}. Both approaches combine candidate forecasts (or estimates) obtained from different sources. Empirical evidence frequently shows the superiority of the combined forecast over the single best forecast for various reasons. For example, combination aggregates the incomplete information \citep{Timmermann2006} and at the same time averages out the error of each candidate forecast caused for instance, by time instability in the specification of models \citep{Rossi2021}. The shrinkage property of combination could also potentially improve forecasting accuracy \citep{Hendry2004} (see \citet{Timmermann2006} for an extensive review.) The literature has witnessed a large and yet increasing number of studies on how to best combine multiple forecasts. Numerous efforts have been devoted to developing data-driven weights in the hope of achieving certain optimality of the combination estimator \citep[see][for a partial list]{Granger1984Improved,Diebold1988serial,Kolassa2011,Hsiao2014,Montero2020}. The combination technique is also extensively studied in a closely related literature on model averaging, where a number of methods have been proposed to determine the weights, such as Bayesian model averaging \citep[see][for a review]{Steel2020}, Mallows' criterion \citep{hansen:2007}, jackknife averaging \citep{hansen.racine:2012}, Kullback-Leibler distance \citep{Zhang2015Kullback}, penalized least squares \citep{zhang2019Parsimonious}, among many others.

	Given candidate forecasts, data-driven combination not only requires researchers to specify which criterion (objective function) to estimate the weights but also in which space one searches for the optimal weight, in other words, which weight constraints should be imposed. A significant portion of the literature has been devoted to answering the first question; see \citet{Wang2023} for an excellent review on this aspect, including the history and recent developments. In contrast, the specification of weight space receives significantly much less attention, and the role of weight constraints on the properties of combination is also less understood, leading to rather arbitrary use of weight constraints in practice. This study offers the first comprehensive review on the weight constraints. We theoretically discuss how various constraints influence the properties of the combined forecast and verify our theory via numerical studies.

	In practice, often-used weight constraints include non-negativity, sum-up-to-unity, norm constraints (see Section~\ref{sec:setup} for precise definitions), among others. Existing studies on forecast combination and model averaging typically employ a (sub)set of these constraints. For example, \cite{Ando2014A,Ando2017A} impose the non-negativity constraint to determine model-averaging weights for high-dimensional models. \citet{li:kang:li:2023} proposes time-varying weighting based on a variant of softmax function which implicitly requires non-negativity. The sum-up-to-unity constraint is advocated by \citet{Diebold1988serial} to eliminate serial correlation in regression-based approaches, and this constraint is also used with the hope of achieving unbiased combination when all candidate forecasts are unbiased \citep[see also, e.g.,][]{Granger1984Improved}.
 Most studies employ the non-negativity and sum-up-to-unity constraints jointly, such as the default weight space in optimal model averaging\citep{hansen:2007, Zhang2016Optimal,CHEN2023,zou2024analyzing,model2025tzu},  smoothed information criteria \citep{hjort.claeskens:2003,claeskens.croux.ea:2006,Rigollet2012sparse}, and averaging based on historical performance, for example, variance and mean squared error. Finally, the norm constraint is often used if the objective function is based on eigenvectors of combined forecasts \citep[see, e.g.,][]{Hsiao2014}.

	Despite its importance in weight estimation, the choice of constraints is far less discussed in the literature. A notable exception is \citet{Radchenko2023} which discusses how the non-negativity constraint plays a role in the combination. Nevertheless, it generally remains unclear to practitioners how the use of individual or multiple of these constraints influences the properties of the combined forecast. Specifically, how does the bias, variance, in-sample and out-of-sample fit of the combined forecast behave when applying different sets of weight constraints? Does a constraint lead to a unique estimated weight? Is the resulting weight sparse, such that only a small number of candidate forecasts eventually contribute to the combination? Lack of a good understanding of these questions leaves the unconscious and perhaps arbitrary choice of weight constraints in practice, further leading to unjustified performance of the combined forecast. This study addresses these questions by theoretically comparing various weight constraints and studying the impact of a set of constraints on the performance of combined forecasts. Inevitably, the impact of weight constraints on the resulting combination forecast is intertwined with weight choice criteria. To facilitate the analysis, we consider several most popular forecast combined methods, including regression-based weights, model-averaging-based weights, performance-based weights, and the eigenvector approaches. We discuss each set of weight constraints paired with every possible \emph{compatible} criterion. Our analysis provides guidance for practitioners to decide which set of weight constraints to use depending on the target.

The rest of this paper is organized as follows. Section \ref{sec:setup} summarizes popular weight constraints used in forecast combination and model averaging.
Section \ref{sec:Weightobtained} presents widely used objective functions for weight estimation in conjunction with constraints.
Section \ref{sec:properties}
analyzes the properties of combined forecasts under different constraints.
Section \ref{sec:guidance} describes two practical ways to determine a proper weight constraint.
Simulation results are provided in Section \ref{sec:simulation}. Finally,
Section \ref{sec:conclusion} concludes this overview with some brief discussion. Proofs are provided in the Appendix.

\section{Forecast combination and weight constraints}\label{sec:setup}
Suppose that we observe $\{y_t,t=1,\ldots,T\}$, and wish to forecast the future values of $y_{T+1}$ by combining $S$ candidate forecasts produced by different models or experts. Let $f_{t,s}$ be the $s$-th candidate forecast at time $t$ for $t\geq T+1$. Denote $\f_t=(f_{t,1},\ldots,f_{t,S})\tp$ as the vector of all candidate forecasts at time $t$, and $\f_{(s)}=(f_{1,s},\ldots,f_{T,s})\tp$ as the $s$-th candidate forecasts for all time horizons.
The final forecast is obtained by combining $\{f_{t,s}\}_{s=1}^S$, that is, $\hat{y}_t=\f_t\tp\w$, where $\w=(w_1,\ldots,w_S)\tp$ is an $S\times 1$ vector of weights.

	The literature has witnessed diversified choices of combination weights with distinct constraints. Here we provide a list of popular weight constraints, under which the weights are optimized. We emphasize that this list is not comprehensive but focuses on the widely used methods in practice that can be analytically analyzed. The benchmark would be no constraints, and we denote this weight space as $\bbW^\calA=\{\w|\w\in\calR^S\}$, which may lead to arbitrarily large weights. To avoid extreme weight values and achieve certain desired statistical properties, a set of weight constraints are typically imposed in practice. First, one can force the weights to sum up to unity, and we denote this weight space as $\bbW^\calB=\{\w| \w\tp\one=1\}$. When each candidate forecast is unbiased, this constraint guarantees the unbiasedness of the combination forecast. It also introduces internal competition among candidate forecasts and alleviates the serial correlation (see Remark \ref{rem:serial}).
 Another widely imposed constraint is non-negativity, that is, $\bbW^\calC=\{\w|\w\in[0,1]^S\}$, making weights more alike probabilities.
 The underlying assumption of constraining weights in the space of $\bbW^\calC$ is that each candidate forecast provides useful information and contributes positively to the final forecast. Combining both sum-up-to-unity and non-negativity constraints, we denote $\bbW^\calD=\left\{\w\big| \w\in [0,1]^S \text{ and } \one\tp\w=1\right\}$. Finally, one can impose a constraint on the norm of weights, namely $\bbW^\calE=\left\{\w\big| \w\tp\w=1\right\}$. This constraint is typically used when combination weights are from an eigenvector-based objective function. Compared with the sum-up-to-unity constraint that restricts the search on a $\calR^{S-1}$ hyperplane, the norm constraint in $\bbW^\calE$ allows the search of the entire $\calR^{S}$ \citep{Hsiao2014}.

\begin{figure}[htp]
\begin{center}
\includegraphics[width=0.45\textwidth]{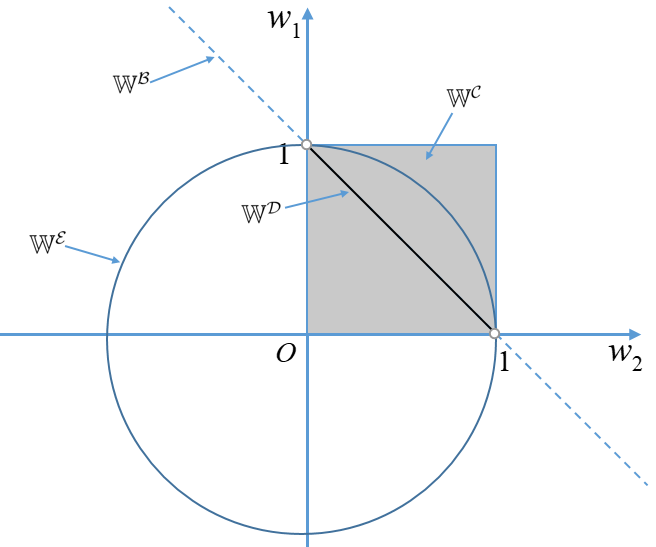}
\end{center}
\caption{The schematic diagram for weight spaces.}\label{fig:regions}
\footnotesize{\emph{Notes:} The sum-up-to-unity weight constraint $\bbW^\calB$ corresponds to the downward sloping 45 degree line passing $(0,1)$ and $(1,0)$. The non-negativity constraint $\bbW^\calC$ restricts weights to be in the shadow box in the upper-right quadrant. Combining sum-up-to-unity and non-negativity constraints $\bbW^\calD$ limits the weights within the dark solid part of the downward sloping 45-degree line. The unity norm constraint $\bbW^\calE$ corresponds to the unit circle.}
\end{figure}

	We can illustrate these four weight constraints via a schematic diagram in a 2-dimensional case (with two candidate forecasts) as Figure~\ref{fig:regions}. The sum-up-to-unity weight space $\bbW^\calB$ corresponds to the downward sloping 45-degree line passing $(0,1)$ and $(1,0)$. The non-negativity constraint $\bbW^\calC$ restricts weights to be in the shadow box in the upper-right quadrant. Combining sum-up-to-unity and non-negativity constraints $\bbW^\calD$ limits the weights within the dark solid part of the downward sloping 45-degree line. Finally, the unity norm constraint $\bbW^\calE$ corresponds to the unit circle.

Many practically popular weight choices fall into the above mentioned weight constraints. For example, the classic forecast combination method by \cite{Granger1984Improved} imposes no constraints and weights are freely chosen from $\bbW^\calA$. The sum-up-to-unity constraint $\bbW^\calB$ is used to eliminate serial correlation in the combination; see, for example, \citet{Diebold1988serial,Diebold1996Forecast,Breiman1996Stacked,Zhou2012Ensemble}.
\cite{Ando2014A,Ando2017A} employ the non-negativity constraint $\bbW^\calC$ to control model-averaging weights for high-dimensional data. The majority of combination methods confine weights to $\bbW^\calD$. These include the simple average, that is, $w_s=1/S$ for $s=1, 2, \ldots,S$ \citep{Clemen1989Combining,Chan1999A}; inverse error weights by \cite{Bates1969}, that is, $w_s=\hat{\sigma}^{-2}_s/\sum_{s=1}^{S}\hat{\sigma}^{-2}_s$,
where $\hat{\sigma}^2_s=T^{-1}\sum_{t=1}^{T}(y_t-f_{t,s})^2$ denoting the estimated mean squared prediction error of the $s$-th candidate model;
and smoothed information criteria \cite[IC, see, e.g.,][]{Hurvich1989Regression,hjort.claeskens:2003,claeskens.croux.ea:2006,Zhang2016Optimal}
\begin{align*}
w_s=\frac{\exp(-\XIC_s/2)}{\sums\exp(-\XIC_s/2)},\quad s=1, \cdots, S,
\end{align*}
where $\XIC_s$ represents a certain IC. Many model averaging methods also restrict weights to $\bbW^\calD$, for example, Mallows averaging \citep{hansen:2007,Fang2022An,chi2025model}, jacknife averaging \citep{hansen.racine:2012,lu2015jackknife} and cross-validation (CV) model averaging \citep{zhang2023model,bu2025improving}. Finally, the norm constraint is adopted by \cite{Hsiao2014} in an eigenvector approach of forecast combination.
	
\begin{remark}
One can also link the non-negativity and sum-up-to-unity constraints with the shrinkage estimator of the covariance matrix of candidate forecasts. We illustrate this link in a simple case where the candidate forecasts are all unbiased and the weights are treated as nonrandom. We aim to minimize the combination variance, that is, $\min_\w \w\tp \bSigma\w$, where $\bSigma$ is the covariance matrix of candidate forecasts. \citet{Radchenko2023} show that when candidate forecasts are highly correlated, the resulting weights without imposing any bound constraint are likely to be negative.
If we impose both sum-up-to-unity and non-negativity constraint, namely $\bbW^\calD$, Proposition 1 of \citet{jagannathan2003risk} implies that a constrained optimum based on $\bSigma$ is equivalent to an unconstrained one obtained from using $\tilde{\bSigma}=\bSigma-(\one\tp \brho+\brho\one\tp)$, where $\brho=(\rho_1,...,\rho_S)\tp$ is the multiplier for the non-negativity constraint.
For the $i$-th candidate forecast, the non-negativity constraint implies that $\Sigma_{is}$ for $s\neq i$ is reduced by $\rho_i+\rho_s$ (a positive quantity), and its variance is reduced by $2\rho_i$. In this sense, the new covariance matrix estimates $\tilde{\bSigma}$ can be regarded as a shrinkage counterpart of the original covariance $\bSigma$.
\end{remark}

\section{Objective functions to determine weights}
\label{sec:Weightobtained}
Admittedly, it is highly difficult, if not impossible to isolate the discussion of weight constraints from the objective function for estimating the weights. Even for the same weight space, different weight estimation methods can lead to substantially different results. However, it is beyond the scope of this paper to review all possible forecast combinations or model averaging methods. Our focus is to compare the effect of various constraints, and thus, we discuss several mostly widely used methods to determine the combination/averaging weights, based on which the weight constraints are imposed. For convenience, we define $\y=(y_1,\ldots,y_T)\tp$ and $\F=(\f_1,\ldots,\f_T)\tp=(\f_{(1)},\ldots,\f_{(S)})$.

\subsection{Regression-based method}
A straightforward method to determine the combination weights is to regress $y_t$ on all candidate forecasts \citep{Granger1984Improved}, namely
$y_t=\f_t\tp \w+\epsilon_t,$
where $\epsilon_t$ is an independently distributed error term with mean zero and variance $\sigma^2$. Then the weight vector can be obtained by $\hat{\w}^\calA_{\reg}=(\F\tp\F)^{-1}\F\tp\y$, where the subscript denotes the estimation method and the superscript presents the weight space. This method is referred to as method A in \cite{Granger1984Improved}.

To guarantee empirical unbiasedness, that is, $\one\tp(\y-\hat{\y})=0$ with $\hat{\y}$ being a forecast of $\y$, \cite{Granger1984Improved} propose to add an intercept in the regression model as
\begin{align}\label{eq:regression-model}
y_t=\delta_0+\f_t\tp\w+\epsilon_t, \quad t=1,\ldots, T.
\end{align}
The resulting weight estimator (called method C in \cite{Granger1984Improved}) is then denoted as $\hat{\w}^{\calA'}_{\reg}$, which can be associated with $\hat{\w}^\calA_{\reg}$ as
\begin{align}\label{eq:wregA2}
\hat{\w}^{\calA'}_{\reg}&=\hat{\w}^\calA_{\reg}-\hat{\delta_0}(\F\tp\F)^{-1}\F\tp\one,
\end{align}
where $\hat{\delta_0}=\theta^{-1}\one\tp\hat{\e}$,
$\theta=n-\one\tp\F(\F\tp\F)^{-1}\F\tp\one$ and $\hat{\e}=\y-\F\hat{\w}^\calA_{\reg}$.

The regression method can also be used jointly with alternative weight constraints. If one considers using the weight space $\bbW^\calB$ in the regression model~\eqref{eq:regression-model}, the resulting weight (called method B in \cite{Granger1984Improved})
can be written as
\begin{align}\label{eq:wregB}
\hat{\w}^\calB_\reg
=\hat{\w}^\calA_{\reg}-\hat{\rho}_0(\F\tp\F)^{-1}\one,
\end{align}
where $ \hat{\rho}_0=(\one\tp \hat{\w}^\calA_{\reg}-1)/\one\tp(\F\tp\F)^{-1}\one$, and this weight helps alleviate the serial correlation (see Remark \ref{rem:serial}). Of course, one can also estimate the weights from regression model~\eqref{eq:regression-model} under the constraint sets $\bbW^\calC$ and $\bbW^\calD$, and obtain $\hat{\w}^\calC_\reg$ and $\hat{\w}^\calD_\reg$, respectively. Unfortunately, these estimates do not have a closed-form solution. For the weight space $\bbW^\calE$, via Lagrangian multiplier method, the optimal solution is $\hat{\w}^\calE_\reg=(\F\tp\F+\hat{\nu})^{-1}\F\tp\y$, where $\hat{\nu}$ satisfies $\y\tp\F(\F\tp\F+\hat{\nu})^{-2}\F\tp\y=1$.

\begin{remark}\label{rem:serial}
\cite{Diebold1988serial} shows that the unrestricted
ordinary least squares (OLS) estimator of regression-based weights introduces serially correlated residuals even if the candidate forecasts have serially uncorrelated
errors \citep[see also][]{de2000review}. To see this, consider the regression model without an intercept, and the combined forecast is given by $\hat{y}_{t}=\f_t\tp\hat{\w}$, where $\hat{\w}$ is a least squares estimator of $\w$, such that the forecast error is
\begin{align}\label{eq:forecast-error}
\hat{y}_t-y_t&=y_t\left(\sums \hat{w}_s-1\right)+\sums \hat{w}_s(f_{t,s}-y_t)\notag\\
&=y_t\left(\sums \hat{w}_s-1\right)+\sums \hat{w}_s\epsilon_{t,s},
\end{align}
where $\epsilon_{t,s}=f_{t,s}-y_t$.
Equation~\eqref{eq:forecast-error} suggests that, if $y_t$ exhibits serial correlation, then the error of the combined forecast is generally serially correlated. Constraining the sum of weights to one alleviates the serial correlation of the combination error.
\end{remark}

\subsection{Model averaging-based method}
Recently, model averaging methods have received increasing attention in dealing with model uncertainty, for example, which regressors to include in a regression model, and various criteria have been proposed to determine the averaging weights. It is conceptually closely related to forecast combination, and one can use model averaging criteria to determine the weights for forecast combination by formulating a regression model of $\y$ on $\F$. We focus on asymptotic optimal model averaging here, because it has a similar goal as forecasting, namely to achieve the best prediction performance.

A prevalent optimal averaging approach is Mallows model averaging \citep{hansen:2007,Fang2022An,chi2025model}. It determines the weight using the Mallows' criterion, which is an unbiased estimator of the risk (ignoring terms that do not depend on weights), but it only works for linear models. Let the $s$-th candidate forecast be a linear projection of the dependent variable, i.e., $\f_\s=\bP_s\y$ with $\bP_s$ being the projection matrix of the $s$-th candidate model. The Mallows criterion can be written as
\begin{align}
\calC(\w)&=\|\y-\F\w\|^2+2\hat{\sigma}^2 \tr\{\bP(\w)\}\notag\\
&=\w\tp\F\tp\F \w+2\w\tp(\hat{\sigma}^2 \bk-\F\tp\y)+\y\tp\y ,\label{eq:formG1}
\end{align}
where $\bk=(\tr(\bP_1),\ldots,\tr(\bP_M))\tp$, $\hat{\sigma}^2$ is the variance of error to approximate $\y$ with $\F\w$, which can be estimated using the full model with all candidate forecasts \citep[see, e.g.,][]{hansen:2007}.

Under linear models, \cite{Zhang2015Kullback} proposes another unbiased optimal averaging criterion based on Kullback-Leibler (KL) divergence,
\begin{align}\label{eq:formG2}
\textrm{KL}(\w)&=\|\y-\F\w\|^2+2\hat{\sigma}^2\tr\{\bP(\w)\}-2\y\tp \bP\tp(\w)\frac{\partial \hat{\sigma}^2}{\partial \y}\notag\\
&=\w\tp\F\tp\F \w+2\w\tp(\hat{\sigma}^2 \bk-\bphi-\F\tp \y)+\y\tp\y,
\end{align}
where $\bk=(\tr(\bP_1),\ldots,\tr(\bP_M))\tp$ and $\bphi=(\y\tp\bP_1 \frac{\partial\hat{\sigma}^2}{\partial\y},\ldots,\y\tp\bP_S \frac{\partial\hat{\sigma}^2}{\partial\y})\tp=(\f_{(1)},\ldots,\f_{(S)})\tp\frac{\partial\hat{\sigma}^2}{\partial\y}$. We can encompass both Mallows and KL criteria in a general framework as
\begin{align}
\calD(\w)=\w\tp\F\tp\F \w+2\w\tp\bpsi+\y\tp\y,\label{eq:objectiveG}
\end{align}
where $\bpsi=\hat{\sigma}^2\bk -\F\tp\y$ for \eqref{eq:formG1} and $\bpsi=\hat{\sigma}^2\bk-\bphi-\F\tp\y$ for \eqref{eq:formG2}. We refer to this criterion as generalized Mallows.

If we impose no restrictions on the weights, namely $\w\in\bbW^\calA$,
then the optimal weight vector can be obtained as
$
\hat{\w}^\calA_{\ma}=-(\F\tp\F)^{-1}\bpsi$.
When the weight is restricted to be in $\bbW^\calB$, solving~\eqref{eq:objectiveG} gives the optimal weight vector as
$
\hat{\w}^\calB_{\ma}=-(\F\tp\F)^{-1}(\bpsi+\check{\rho}_0\one),
$
where
$\check{\rho}_0=-\{\bpsi\tp(\F\tp\F)^{-1}\one+1\}/\one\tp (\F\tp\F)^{-1}\one$.
When the weight belongs to $\bbW^\calC$, there is generally not a closed-form solution, because we cannot determine which boundary condition is binding, but we denote that optimal weight as $\hat{\w}^\calC_{\ma}$.
Imposing weight constraints $\bbW^\calD$ on the averaging criterion in \eqref{eq:objectiveG} produces the weight $\hat{\w}^\calD_{\ma}$ that also lacks a closed-form.
One can also impose the constraint $\bbW^\calE$ to~\eqref{eq:objectiveG}. By the Lagrangian multiplier method, the optimal weight can be obtained by
$
\hat{\w}^\calE_{\ma}=-(\F\tp\F+\check{\nu})^{-1}\bpsi,
$
where $\check{\nu}$ satisfies $\bpsi\tp(\F\tp\F+\check{\nu})^{-2}\bpsi=1$.

Alternative to Mallows or KL criterion, if one is ignorant about the distribution of data, the cross-validation or jackknife method is often used to determine the optimal averaging weights \citep[see, e.g.,][]{hansen.racine:2012,Zhang2013jackknife,lu2015jackknife,Zhang2020Cross}. The leave-one-out cross-validation (CV) criterion minimizes the following objective function:
\begin{align}\label{eq:formH}
\CV(\w)=\sum_{t=1}^T(y_t-\f^{[-t]\top}_t\w)^2=\|\y-\bar{\F}\w\|^2,
\end{align}
where $\f_t^{[-t]}=(f_{t,1}^{[-t]},\ldots,f_{t,S}^{[-t]})\tp$ is the vector of candidate forecasts without using the $i$-th observation, and $\bar{\F}=(\f_1^{[-1]},\ldots,\f_T^{[-T]})\tp$.
Imposing no constraints, we can obtain weights from \eqref{eq:formH} as
$\hat{\w}^\calA_{\cv}=(\bar{\F}\tp\bar{\F})^{-1}\bar{\F}\tp\y.
$ When we impose the constraint $\bbW^\calB$,
the resulting weight is:
$\hat{\w}^\calB_\cv
=\hat{\w}^\calA_{\cv}-\bar{\rho}_0(\bar{\F}\tp\bar{\F})^{-1}\one$,
where $\bar{\rho}_0=(\one\tp \hat{\w}^\calA_{\cv}-1)/\one\tp(\bar{\F}\tp\bar{\F})^{-1}\one$.
When the weight belongs to $\bbW^\calC$ or $\bbW^\calD$, again the resulting weights $\hat{\w}^\calC_\cv$ and $\hat{\w}^\calD_\cv$ do not have a closed-form solution. Under the constraint $\bbW^\calE$, we can obtain the weight from~\eqref{eq:formH} as
$
\hat{\w}^\calE_\cv=(\bar{\F}\tp\bar{\F}+\bar{\nu})^{-1}\bar{\F}\tp\y
$,
where $\bar{\nu}$ satisfies $\y\tp\bar{\F}(\bar{\F}\tp\bar{\F}+\bar{\nu})^{-2}\bar{\F}\tp\y=1$.

\subsection{Individual performance-based methods}
The individual performance-based method typically aims to achieve the best performance by combining forecasts based on a certain measure of their historical performance. \cite{zhang2010phd} proposes a general form of individual performance-based weights,  namely
\begin{align}\label{eq:25}
w_s=\frac{a^{q_s}(n-q_s)^b(\hat{\sigma}_s^2)^c}{\sumjs a^{q_j}(n-q_j)^b(\hat{\sigma}_j^2)^c},
\end{align}
where $a>0, b\geq 0, c\leq 0$, $q_j\geq 0$ and $\hat{\sigma}_s^2$ is the maximum likelihood estimator
of the variance of the $s$-th candidate forecast.
When $a=e^{-1}, b=0$ and $c=-n/2$, \eqref{eq:25} gives the smoothed AIC weights and when $a=n^{-1/2}, b=0$ and $c=-n/2$, it reduces to the smoothed BIC weights, both of which take the form
$\exp(-\textrm{IC}_s/2)/\sumjs\exp(-\textrm{IC}_j/2)$, with IC$_s$ being either AIC or BIC, for the $s$-th forecast \citep{Buckland:1997}.

Besides, one can design a weighting method based on the inverse of a certain loss function, such that the weight of the $s$-th candidate forecast takes a general form as
\begin{equation}\label{eq:inverse-performance}
w_s=\frac{L_s^{-1}}{\sumjs L_j^{-1}},
\end{equation}
where $L_s$ is a loss function of the $s$-th forecast. For example, \citet{Bates1969} measures the performance via mean residual sum, defined as $\hat{\sigma}^2_s n/(n-q_s)$ for the $s$-th candidate forecast, and the weight obtained under this measure can be viewed as a special case of~\eqref{eq:25} when $a=b=1$ and $c=-1$.
\cite{Stock1998AComparison} considers individual performance-based weights based on
mean squared error (MSE) in a rolling window manner, which can be written by using the MSE as the loss function in~\eqref{eq:inverse-performance}. \cite{NOWOTARSKI2014Anempirical} measures the performance via the root mean squared errors (RMSE), while \cite{aiolfi2006persistence} and \cite{andrawis2011combination} consider the performance rank.

Obviously, these individual performance-based methods all constrain the weights to be in the space $\bbW^\calD$, and thus we denote this category of weights as $\hat{\w}^\calD_{\pf}$

\subsection{Eigenvector approach}
\cite{Hsiao2014} introduces an eigenvector approach that determines the combination weight by:
\begin{align*}
\min_{\w}\frac{1}{T}\sum_{t=1}^{T}\left\{(y_t\one-\f_t)\tp\w\right\}^2(=T^{-1}\left\|(\y\otimes \one\tp-\F)\w\right\|^2),\quad\text{s.t.}\quad \w\in\bbW^\calE,
\end{align*}
where $\otimes$ is the Kronecker product. The resulting weight vector, denoted as $\hat{\w}^{\calE}_{\eig}$, is the eigenvector belonging to the smallest eigenvalue of $\M=T^{-1}\sum_{t=1}^{T} (y_t\one-\f_t)(y_t\one-\f_t)\tp=T^{-1}(\y\otimes \one\tp-\F)\tp(\y\otimes \one\tp-\F)$.
The main motivation of the eigenvector approach is to treat the uncertainties in $y_t$ and $\f_t$ symmetrically by attaching weights to the forecast error, aiming to
achieve the geometrically ``best'' fit of the subspace
to the points $y_t\one-\f_t$ for all $t$. This is in sharp contrast to the regression-based method that only attaches weights to $\f_t$, implicitly assuming that there is no uncertainty in $\f_t$ but only in $y_t$. Thus, \citet{Hsiao2014} argues that the eigenvector approach is expected to be less sensitive to the outlying observations of $y_t$, and the resulting weights are also less likely to take extremely large values. Compared with $\bbW^\calD$, the resulting weight from $\bbW^\calE$ is usually not sparse, such that many candidate models can contribute to the combination.

\section{Properties of weight constraints}\label{sec:properties}
This section examines the properties of different weight constraints. These properties unavoidably depend on the weight estimation methods. Thus, to facilitate analysis, we compare the constraints under each category of estimation methods. Of course, some constraints are only relevant for certain estimation methods, for example, individual performance-based method implies $\bbW^\calD$.

\subsection{The sum of squared residuals}\label{sec:SSR}
Following \cite{Granger1984Improved}, we intend to compare the different weight constraints on the fitness of training data. One of the most common measures of fitness is the sum of squared residuals (SSR), defined as $\|\y-\hat{\y}\|^2=\sum_{t=1}^{T}(y_t-\f_t\tp\hat{\w})^2$.

We first consider regression-based methods and examine the weights with an analytical form, namely $\hat{\w}^{\calA}_\reg$, $\hat{\w}^{\calA\prime}_\reg$ and $\hat{\w}^\calB_\reg$.
Based on \eqref{eq:wregA2} and \eqref{eq:wregB}, the SSRs of these weights can be obtained as
\begin{align*}
\SSR^\calA_\reg
&=\|\y-\F\hat{\w}^\calA_\reg\|^2
=\y\tp\{\I-\F(\F\tp\F)^{-1}\F\tp\}\y,\\
\SSR^{\calA'}_\reg
&=\|\y-\hat{\delta}_0\one-\F\hat{\w}^{\calA'}_\reg\|^2\notag\\
&=\|\y-\F\hat{\w}^\calA_\reg\|^2-2\hat{\delta}_0 \one\tp\left\{\I-\F(\F\tp\F)^{-1}\F\tp\right\}\hat{\e}\notag\\
&\quad+\hat{\delta}_0^2\one\tp \left\{\I-\F(\F\tp\F)^{-1}\F\tp\right\} \left\{\I-\F(\F\tp\F)^{-1}\F\tp\right\}\one \notag\\
&=\SSR^\calA_\reg-2\hat{\delta}_0\one\tp
\{\I-\F(\F\tp\F)^{-1}\F\tp\}\y+\hat{\delta}_0^2
\one\tp\{\I-\F(\F\tp\F)^{-1}\F\tp\}\one\notag\\
&=\SSR^\calA_\reg-\theta\hat{\delta}_0^2,
\end{align*}
and
\begin{align}\label{eq:SSRB}
\SSR^\calB_\reg&=\|\y-\F\hat{\w}^\calB_\reg\|^2\notag\\
&=\|\y-\F\hat{\w}^\calA_\reg\|^2+2(\y-\F\hat{\w}^\calA_\reg)\tp\{\hat{\rho}_0\F(\F\tp\F)^{-1}\one\}+\|\hat{\rho}_0\F(\F\tp\F)^{-1}\one\|^2\notag\\
&=\SSR^\calA_\reg+\hat{\rho}_0^2\{\one\tp(\F\tp\F)^{-1}\one\},
\end{align}
where
$\hat{\delta}_0=\theta^{-1}\one\tp\hat{\e}$,
$\theta=n-\one\tp\F(\F\tp\F)^{-1}\F\tp\one$,
$\hat{\e}=\y-\F(\F\tp\F)^{-1}\F\tp\y$ and $\hat{\rho}_0=(\one\tp \hat{\w}^\calA_\reg-1)/\one\tp(\F\tp\F)^{-1}\one$.
Comparing these three SSRs, we can show that $$
\SSR^\calB_\reg\geq \SSR^\calA_\reg\geq \SSR^{\calA'}_\reg.$$
For the weights that lack a closed-form solution, we cannot derive the resulting SSRs explicitly. However, we can infer their relationship by examining the respective optimization problem. Note that in general for a given objective function, $f(\w)$, we have $\min_{\w\in\bbA}f(\w) \geq \min_{\w\in\bbB} f(\w)$ if the solution space $\bbA\subset \bbB$.
Since regression-based methods directly target minimizing the SSR and the weight space satisfies the following relation: $\bbW^\calD\subset\bbW^\calC\subset \bbW^\calA$ , $\bbW^\calD\subset\bbW^\calB\subset \bbW^\calA$ and $\bbW^\calE\subset\bbW^\calA$, we can obtain the following relations:
$$
\SSR^\calD_\reg\geq \SSR^\calC_\reg\geq\SSR^\calA_\reg,\quad
\SSR^\calD_\reg\geq\SSR^\calB_\reg\geq\SSR^\calA_\reg,\quad
\SSR^\calE_\reg\geq\SSR^\calA_\reg.
$$

Next, we consider the optimal averaging method. 	While the Mallows criterion in \eqref{eq:formG1} does not directly target the SSR, its expectation (asymptotically) equals the expectation of the regression-based method if we ignore the term that does not depend on $\w$ \citep{hansen:2007}. Thus, the SSR comparison across weight constraints using optimal averaging remains similar to regression-based methods, that is,
$$
\SSR^\calD_\ma\geq \SSR^\calC_\ma\geq\SSR^\calA_\ma,\quad
\SSR^\calD_\ma\geq\SSR^\calB_\ma\geq\SSR^\calA_\ma,\quad
\SSR^\calE_\ma\geq \SSR^\calA_\ma.
$$
From a different perspective, we can also compare the SSR of optimal averaging weights via the link with the regression-based weights. We can show that:
\begin{align}\label{eq:SSRAma}
\SSR^\calX_\ma
&=\|\y-\F\hat{\w}^\calX_\ma\|^2\notag\\
&=\|\y-\F\hat{\w}^\calA_\ma\|^2-2\y\tp\F(\hat{\w}^\calX_\ma-\hat{\w}^\calA_\reg)
+
2\hat{\w }^{\calA\top}_\reg  \F\tp\F(\hat{\w}^\calX_\ma-\hat{\w}^\calA_\reg)\notag\\
&\quad+
(\hat{\w}^\calX_\ma-\hat{\w}^\calA_\reg)\tp\F\tp\F(\hat{\w}^\calX_\ma-\hat{\w}^\calA_\reg)\notag\\
&=\SSR^\calA_\reg+(\hat{\w}^\calX_\ma-\hat{\w}^\calA_\reg)\tp\F\tp\F(\hat{\w}^\calX_\ma-\hat{\w}^\calA_\reg),
\end{align}
where $\calX=\calA, \ldots, \calE$. From \eqref{eq:SSRAma}, we have $\SSR_\ma^\calX\leq \SSR_\ma^\calY$ if $\|\hat{\w}^\calX_\ma-\hat{\w}^\calA_\reg\|\leq \|\hat{\w}^\calY_\ma-\hat{\w}^\calA_\reg\|$ for $\calX, \calY=\calA, \ldots, \calE$. This result suggests that the closer a Mallows averaging weight to $\hat{\w}^\calA_\reg$, the smaller SSR it produces.

For CV model averaging in \eqref{eq:formH}, following \eqref{eq:SSRB}, we have:
\begin{align*}
\SSR^\calB_\cv&=\|\y-\F\hat{\w}^\calB_\cv\|^2\notag\\
&=\|\y-\F\hat{\w}^\calA_{\cv}\|^2+2(\y-\F\hat{\w}^\calA_{\cv})\tp\{\bar{\rho}_0\F(\bar{\F}\tp\bar{\F})^{-1}\one\}+\|\bar{\rho}_0\F(\bar{\F}\tp\bar{\F})^{-1}\one\|^2\notag\\
&\approx\SSR^\calA_\cv+\bar{\rho}_0^2\one\tp(\bar{\F}\tp\bar{\F})^{-1}\one,
\end{align*}
where the last equality is due to $\bar{\F}\tp\bar{\F}\approx {\F}\tp{\F}$ because the omitted ones become negligible when $T$ is large.
This implies that $\SSR^\calB_\cv\geq \SSR^\calA_\cv$.
Since the CV objective function is also a sum of squared residuals (but with leave-one-out candidate forecasts), using similar arguments as regression-based methods, we have the following relation:
\begin{align}
\SSR^\calD_\cv\geq \SSR^\calC_\cv\geq\SSR^\calA_\cv, \quad
\SSR^\calD_\cv\geq\SSR^\calB_\cv\geq\SSR^\calA_\cv,\quad
\SSR^\calE_\cv\geq\SSR^\calA_\cv.
\end{align}

Finally, since the individual performance-based method and the eigenvector approach both imply a specific weight space, namely $\hat{\w}^\calD_\pf\in\bbW^\calD$ and $\hat{\w}^{\calE}_\eig\in\bbW^\calE$, we do not compare different weight constraints for these two methods.

Note that the regression-based method directly minimizes the SSR objective function, and thus, given the same weight space, it is expected to produce the minimum SSR than other methods that do not target the SSR, such as the individual performance-based methods, (generalized) Mallows model averaging, and the eigenvector approach.

\subsection{Empirical unbiasedness}
Another important criterion to evaluate the fitness of training data is the empirical unbiasedness \citep{Granger1984Improved}, defined as $\one\tp(\y-\hat{\y})=0$. If the training data are randomly generated from a common distribution, the empirical unbiasedness also implies the asymptotic unbiasedness.

In practice, the empirical unbiasedness is difficult to achieve unless two sufficient conditions are satisfied: (1) the error of each candidate forecast has a zero mean, that is, $\one\tp\y=\one\tp\f_{(s)}$ for $s=1,\ldots,S$; (2) the combination weights add up to unity, that is, $\one\tp \hat{\w}=1$. In this sense, the weights resulting from weight space $\bbW^\calB$ and $\bbW^\calD$ satisfy the second condition, and the combined forecast will be unbiased if each candidate forecast is unbiased. However, the weights obtained from $\bbW^\calA$, $\bbW^\calC$ and $\bbW^\calE$ generally cannot achieve empirical unbiasedness even under unbiased candidate forecasts, except $\hat{\w}^{\calA'}_\reg$ which corrects the bias by including an intercept in the regression model.

\subsection{Conditional mean squared forecasting error}
While the in-sample fit is a relevant criterion, the out-of-sample fit is of more practical interest for the forecasting purpose. A common measure of out-of-sample fit is the mean squared forecasting error (MSFE). We first analyze the MSFE assuming that the candidate forecasts $\{\f_t\}_{t=1}^{T+1}$ are given, such that the randomness solely comes from $\y$. The fixed forecasts assumption can be partially justified by conditioning on the forecasts, and we will relax this assumption in the next subsection.
Denote $\Exp_*(\cdot)$, $\Var_*(\cdot)$ and $\Cov_*(\cdot)$ as the conditional expectation, variance and covariance, respectively, for example, $\Exp_*(\cdot)=\Exp(\cdot|\f_1,\cdots,\f_{T+1})$, then conditional MSFE given candidate forecasts can be written as
\begin{align}
&\quad\Exp_*(y_{T+1}-\hat{y}_{T+1})^2\notag\\
&=\Exp_*\left\{y_{T+1}-\mu_{T+1}+\mu_{T+1}-\Exp_*(\hat{y}_{T+1})+
\Exp_*(\hat{y}_{T+1})-\hat{y}_{T+1}\right\}^2\notag\\
&=\Exp_*(y_{T+1}-\mu_{T+1})^2+\{\mu_{T+1}-\Exp_*(\hat{y}_{T+1})\}^2+
\Exp_*\left\{\Exp_*(\hat{y}_{T+1})-\hat{y}_{T+1}\right\}^2\notag\\
&=\sigma^2+\{\mu_{T+1}-\Exp_*(\hat{y}_{T+1})\}^2+\Var_*(\hat{y}_{T+1}),\label{eq:mse}
\end{align}
where $\mu_{T+1}=\Exp_*(y_{T+1})$.
From the last equality, we see that the conditional MSFE of the combined forecast depends on three terms. The first term $\sigma^2$ is the variance of error disturbance that is common for any combination method. The second term $\{\mu_{T+1}-\Exp_*(\hat{y}_{T+1})\}^2$ measures the squared bias of the combined forecast, and the third term $\Var_*(\hat{y}_{T+1})$ is the conditional forecasting variance. Both the conditional bias and variance (and thus the second and third terms) depend on which combination method is used.

We now examine the bias for different combination methods and weight constraints. Assume that there is a weight vector $\w_{0}$ and a constant $\delta_0$ such that $\Exp_*(y_t)=\delta_0+\f_{t}\tp\w_{0}$ for $t=1,\cdots,T+1$. This assumption is necessary because if $\mu_{t}$ cannot be expressed as a linear combination of $\{f_{t,s}\}_{s=1}^S$ for $t=1, 2, \ldots, T+1$, it implies that a weight to recover the true conditional mean of $y_t$ does not exist, and the difference between the conditional mean of the true value and the combined forecast, namely $\delta_0$, can be arbitrarily complicated, making it difficult to analyze the MSFE.
	
We first examine the bias of different estimation methods and constraints. For the regression-based method, the unconstrained weight $\hat{\w}^\calA_\reg$ produces the bias as
\begin{align}\label{eq:28}
(\bias^{\calA}_\reg)^2=\{\mu_{T+1}-\Exp_*(\hat{y}^{\calA}_{\reg,T+1})\}^2
=\left[\delta_0+\f\tp_{T+1}\{\w_{0}-\Exp_*(\hat{\w}^\calA_\reg)\}\right]^2,
\end{align}
where $\hat{y}^\calA_{\reg, T+1}=\f_{T+1}\tp\hat{\w}^\calA_\reg$. Equation~\eqref{eq:28} shows that the magnitude of bias is mainly determined by the bias of $\hat{\w}^\calA_\reg$, that is, $|\w_{0}-\Exp_*(\hat{\w}^\calA_\reg)|$, in which the conditional expectation of weights $\Exp_*(\hat{\w}^\calA_\reg)$ can be written as
\begin{align}\label{eq:29}
\Exp_*(\hat{\w}^\calA_\reg)&=\Exp_*\left\{(\F\tp\F)^{-1}\F\tp\y\right\}\notag\\
&=(\F\tp\F)^{-1}\F\tp\bmu\notag\\
&=(\F\tp\F)^{-1}\F\tp(\delta_0\one+\F\w_{0})\notag\\
&=\delta_0(\F\tp\F)^{-1}\F\tp\one+\w_{0}.
\end{align}
Combining \eqref{eq:28} and \eqref{eq:29}, the combined forecast using $\hat{\w}^\calA_\reg$ is unbiased if only if $\delta_0=0$. The bias precisely depends on the magnitude of $\delta_0$.
To remove the bias even under nonzero $\delta_0$, one can include an intercept in the regression model for weight estimation, namely using $\hat{\w}^{\calA'}_\reg$. In this case, the conditional bias of the combined forecast can be written as
\begin{align*}
(\bias^{\calA'}_\reg)^2=\{\mu_{T+1}-\Exp_*(\hat{y}^{\calA'}_{\reg,T+1})\}^2
=\left[\delta_0-\Exp_*(\hat{\delta}_0)+\f\tp_{T+1}\{\w_{0}-\Exp_*(\hat{\w}^{\calA'}_\reg)\}\right]^2.
\end{align*}
Due to the unbiasedness of least squares estimation, $\Exp_*(\hat{\delta}_0)=\delta_0$ and $\Exp_*(\hat{\w}^{\calA'}_\reg)=\w_{0}$, and thus we have $\bias^{\calA'}_\reg=0$.

When constraints are imposed, the resulting weights usually do not have a closed-form solution, making it more difficult to analyze their bias. If $\delta_0=0$, the main part of bias is $|\Exp_*(\hat{\w}^{\calX})-\w_{0}|$, where $\calX$ is $\calA, \calB, \calC$, $\calD$ or $\calE$. Since the objective function in the regression-based method is a quadric loss, $\Exp_*(\hat{\w})$ is usually the smallest-distance approximation of $\w_{0}$ in feasible regions, that is,
$|\Exp_*(\hat{\w}^\calX_{\cdot})-\w_{0}|\approx\inf_{\w\in\bbW^\calX}|\w-\w_{0}|$. Further noting that for any two weight spaces $\mathbb{M}$ and $\mathbb{N}$, if $\mathbb{M}\subset\mathbb{N}$ then $\inf_{\w\in \mathbb{M}}|\w-\w_{0}| \geq \inf_{\w\in\mathbb{N}}|\w-\w_{0}|$. Thus, based on the fact that  $\bbW^\calD\subset\bbW^\calB\subset\bbW^\calA$, $\bbW^\calD\subset\bbW^\calC\subset\bbW^\calA$ and $\bbW^\calE\subset\bbW^\calA$, we have:
\begin{gather*}
\bias^\calD_\reg\geq  \{\bias^\calB_\reg, \bias^\calC_\reg\} \geq \bias^\calA_\reg \geq \bias^{\calA'}_\reg\quad\text{and}\quad \bias^\calE_\reg\geq \bias^\calA_\reg.
\end{gather*}
If $\delta_0\neq 0$, the combined forecast is generally biased except $\hat{\w}^{\calA'}_\reg$. Note that we can write the bias as
\begin{align*}
(\bias^{\calX}_\reg)^2&=\{\mu_{T+1}-\Exp_*(\hat{y}^{\calX}_{\reg,T+1})\}^2\notag\\
&=\left[\delta_0+\f\tp_{T+1}\{\w_{0}-\Exp_*(\hat{\w}^\calX_\reg)\}\right]^2\notag\\
&\leq 2\delta_0^2+2\left[\f\tp_{T+1}\{\w_{0}-\Exp_*(\hat{\w}^\calX_\reg)\}\right]^2\notag\\
&\leq 2\delta_0^2+2\|\f_{T+1}\|^2\|\w_{0}-\Exp_*(\hat{\w}^\calX_\reg)\|^2,
\end{align*}
for $\calX = \calA$, $\calB$, $\calC$, $\calD$ or $\calE$. Thus, the presence of $\delta_0$ potentially inflates the upper bound of the bias.

Since the objective function of CV averaging converges to the loss function of the regression-based method, similar bias properties apply to CV averaging. Particularly, when no constraint is imposed, we have:
\begin{align*}
\Exp_*(\hat{\w}^\calA_\cv)&=\Exp_*\left\{(\bar{\F}\tp\bar{\F})^{-1}\bar{\F}\tp\y\right\}\notag\\
&=(\bar{\F}\tp\bar{\F})^{-1}\bar{\F}\tp\bmu\notag\\
&=(\bar{\F}\tp\bar{\F})^{-1}\bar{\F}\tp(\delta_0\one+\F\w_{0})\notag\\
&\approx \delta_0(\bar{\F}\tp\bar{\F})^{-1}\bar{\F}\tp\one+\w_{0},
\end{align*}
where the last equality is due to $\bar{\F}\tp{\F}\approx \bar{\F}\tp\bar{\F}$  because the omitted ones become negligible when $T$ is large;
and the (un)biasedness of $\hat{\w}^\calA_\cv$ depends on $\delta_0$.
When the sum-to-unity constraint is imposed, $\hat{\w}^\calB_\cv$ is also conditionally biased because:
\begin{align*} \Exp_*(\hat{\w}^\calB_\cv)
&=\Exp_*\left\{(\bar{\F}\tp\bar{\F})^{-1}\bar{\F}\tp\y-\bar{\rho}_0(\bar{\F}\tp\bar{\F})^{-1}\one\right\}\notag\\ &=(\bar{\F}\tp\bar{\F})^{-1}\bar{\F}\tp\bmu-\Exp_*(\bar{\rho}_0)(\bar{\F}\tp\bar{\F})^{-1}\one\notag\\ &=(\bar{\F}\tp\bar{\F})^{-1}\bar{\F}\tp(\delta_0\one+\F\w_{0})-\Exp_*(\bar{\rho}_0)(\bar{\F}\tp\bar{\F})^{-1}\one\notag\\ &\approx \w_{0}+\delta_0(\bar{\F}\tp\bar{\F})^{-1}\bar{\F}\tp\one-\Exp_*(\bar{\rho}_0)(\bar{\F}\tp\bar{\F})^{-1}\one,
\end{align*}
where the last equality is due to $\bar{\F}\tp{\F}\approx \bar{\F}\tp\bar{\F}$  because the omitted ones become negligible when $T$ is large.
With similar arguments as in the regression-based method, we can also obtain a similar relation of bias under different weight constraints as regression-based methods, namely
\begin{gather*}
\bias^{\calD}_\cv\geq \{\bias^{\calC}_\cv,\bias^{\calB}_\cv\}\geq \bias^{\calA}_\cv,\quad\text{and}\quad
\bias^{\calE}_\cv\geq \bias^{\calA}_\cv.
\end{gather*}

The generalized Mallows averaging considers a different objective function rather than the quadratic loss, leading to different bias properties. When no constraint is imposed, we have:
\begin{align}
\Exp_*(\hat{\w}^\calA_\ma)
&=-\Exp_*\left\{(\F\tp\F)^{-1}\bpsi\right\}\notag\\
&=\Exp_*\left[(\F\tp\F)^{-1}\{\F\tp\y+\Exp_*(\bphi)-\hat{\sigma}^2\bk\}\right]\notag\\
&=(\F\tp\F)^{-1}\{\F\tp\bmu+\Exp_*(\bphi)-\hat{\sigma}^2\bk\}\notag\\
&=(\F\tp\F)^{-1}\left[\F\tp(\delta_0\one+\F\w_{0})+\Exp_*(\bphi)-\hat{\sigma}^2\bk\right]\notag\\
&=\w_{0}+(\F\tp\F)^{-1}\{\delta_0\F\tp\one-\Exp_*(\bphi)+\hat{\sigma}^2\bk\},\label{eq:66}
\end{align}
where $\bphi=\mathbf{0}$ for Mallows averaging \eqref{eq:formG1} and $\bphi\neq\mathbf{0}$ for KL averaging \eqref{eq:formG2}.
Equation \eqref{eq:66} suggests that $\hat{\w}^\calA_\ma$ is not conditionally unbiased even though $\delta_0=0$.
Under the weight constraints $\bbW^\calB$, we have:
\begin{align}\label{eq:59}
\Exp_*(\hat{\w}^\calB_\ma)&=-\Exp_*\{(\F\tp\F)^{-1}(\bpsi+\check{\rho}_0\one)\}\notag\\
&=-(\F\tp\F)^{-1}\{\Exp_*(\bpsi)+\Exp_*(\check{\rho}_0)\one\}\notag\\
&=\{\Exp_*(\bphi)-\sigma^2\bk\}(\F\tp\F)^{-1}+(\F\tp\F)^{-1}\F\tp\bmu-\Exp_*(\check{\rho}_0)(\F\tp\F)^{-1}\one\notag\\
&=\{\Exp_*(\bphi)-\sigma^2\bk\}(\F\tp\F)^{-1}+(\F\tp\F)^{-1}\F\tp(\delta_0\one+\F\w_{0})-\Exp_*(\check{\rho}_0)(\F\tp\F)^{-1}\one\notag\\
&=\w_{0}+\{\Exp_*(\bphi)-\sigma^2\bk\}(\F\tp\F)^{-1}+\delta_0(\F\tp\F)^{-1}\F\tp\one-\Exp_*(\check{\rho}_0)(\F\tp\F)^{-1}\one,
\end{align}
where $\check{\rho}_0=-\{\bpsi\tp(\F\tp\F)^{-1}\one+1\}/\one\tp (\F\tp\F)^{-1}\one$.
Thus, $\hat{\w}^\calB_\ma$ is also generally conditionally biased because the last three terms of \eqref{eq:59} are nonzero.

Due to the non-quadratic feature of the objective function of generalized Mallows' averaging, it is difficult to associate the weight constraints with bias, even assuming $\Exp_*(y_t)=\f_t\tp\w_0$.
Thus, we analyze the (rough) upper bound of bias as follows. Note that for a general weight $\w$, we have:
\begin{align}
\{\mu_{T+1}-\Exp_*(\hat{y}_{T+1})\}^2
&=\{\mu_{T+1}-\f_{T+1}\tp\Exp_*(\hat{\w})\}^2\notag\\
&\leq 2\mu^2_{T+1}+2\f_{T+1}\tp\Exp_*(\hat{\w})\Exp_*(\hat{\w}\tp)\f_{T+1}\notag\\
&\leq 2\mu^2_{T+1}+2\|\Exp_*(\hat{\w})\|^2\f_{T+1}\tp\f_{T+1}\notag\\
&\leq
\begin{cases}
\infty, & \mbox{if } \hat{\w}\in\bbW^\calA\text{ or }\bbW^\calB \\
2\mu^2_{T+1}+2S\f_{T+1}\tp\f_{T+1}  , & \mbox{if } \hat{\w}\in\bbW^\calC \\
2\mu_{T+1}^2+2\f_{T+1}\tp\f_{T+1}, & \mbox{if } \hat{\w}\in\bbW^\calD\text{ or }\bbW^\calE
\end{cases},
\end{align}
where $S$ is the number of candidate models.
The bound analysis shows that the unconstrained and sum-up-to-unity weight can be biased without an upper bound, while the bound of weight constraint $\bbW^\calD$ and $\bbW^\calE$ is typically smaller than that of $\bbW^\calC$. We summarize the bias of different constraints in the following proposition.
\begin{proposition}\label{pro:bias}\
\begin{itemize}
\item[(1)]
If there exists a weight vector $\w_{0}$ and a $\delta_0\neq 0$ such that $\Exp_*(y_t)=\delta_0+\f_t\tp\w_{0}$ for $t=1, \cdots, T+1$, then the combined forecast is biased except using $\hat{\w}^{\calA'}_\reg$.
\item[(2)]
If there is a weight vector $\w_{0}$ such that $\Exp_*(y_t)=\f_t\tp\w_{0}$ for $t=1, \cdots, T+1$, then
\begin{gather*}
\bias^\calD_\reg\geq \{ \bias^\calB_\reg, \bias^\calC_\reg \}\geq \bias^\calA_\reg \geq \bias^{\calA'}_\reg \quad\text{and}\quad
\bias^\calE_\reg\geq \bias^\calA_\reg;\\
\bias^{\calD}_\cv\geq \{\bias^{\calB}_\cv,\bias^{\calC}_\cv\}\geq \bias^{\calA}_\cv,\quad\text{and}\quad
\bias^{\calE}_\cv\geq \bias^{\calA}_\cv.
\end{gather*}
\item[(3)] The upper bound of conditional bias under different weight spaces is
\begin{align*}
\{\mu_{T+1}-\Exp_*(\hat{y}_{T+1})\}^2\leq
\begin{cases}
\infty, & \mbox{if } \hat{\w}\in\bbW^\calA\text{ or }\bbW^\calB \\
2\mu^2_{T+1}+2S\f_{T+1}\tp\f_{T+1}  , & \mbox{if } \hat{\w}\in\bbW^\calC \\
2\mu_{T+1}^2+2\f_{T+1}\tp\f_{T+1}, & \mbox{if } \hat{\w}\in\bbW^\calD\text{ or }\bbW^\calE
\end{cases}.
\end{align*}
\end{itemize}
\end{proposition}

	Next, we compare the variance of different combined forecasts in \eqref{eq:mse}. As in the bias analysis, we study the exact variance relation if the weights have a closed form, whereas we examine the upper bound of the variance if a closed-form solution of weights is not available. Note that the upper bound of variance is mainly determined by the constraints imposed. A tighter constraint is typically associated with a smaller upper bound of the variance since it limits the variability of estimated weights. We summarize the upper bound of variance resulting from different constraints in the following proposition. 
\begin{proposition}\label{pro:var}\
\begin{itemize}
\item [(1)]$ \Var_*(\hat{y}_{\reg,T+1}^{\calA'})\geq \Var_*(\hat{y}^\calA_{\reg,T+1})\geq \Var_*(\hat{y}_{\reg,T+1}^{\calB})$
and $\Var_*(\hat{y}_{\cv,T+1}^{\calA})\geq \Var_*(\hat{y}_{\cv,T+1}^\calB)$.
\item[(2)] $\Var_*(\hat{y}_{Z,T+1}^\calC)\leq S\f\tp_{T+1}\f_{T+1}$, where $Z$ represents $\reg$, $\ma$ and $\cv$.
\item[(3)]$\Var_*(\hat{y}_{Z,T+1}^\calD)\leq  \f\tp_{T+1}\f_{T+1}$, where $Z$ represents $\reg$,  $\ma$, $\cv$ and $\pf$.
\item[(4)]$\Var_*(\hat{y}_{\eig,T+1}^\calE)\leq  \f\tp_{T+1}\f_{T+1}$.
\end{itemize}
\end{proposition}
Proof. See Appendix \hyperref[sec:appendix_A]{A}.
	
From Propositions \ref{pro:bias} and \ref{pro:var} jointly, we find that a certain type of constraint typically imposes opposite effects on bias and variance. Generally, the combination variance typically increases when fewer (restricted) constraints are imposed and a larger degree of freedom is allowed, which, on the other hand, reduces the bias. This result suggests a typical bias-variance is involved when a weight constraint is imposed. Based on the bound analysis of bias and variance, we can obtain the upper bound of the conditional MSFE as
\begin{align*}
\Exp_*(y_{T+1}-\hat{y}_{T+1})^2
\leq
\begin{cases}
\infty, & \mbox{if } \hat{\w}\in\bbW^\calA\text{ or }\bbW^\calB \\
2\mu^2_{T+1}+3S\f_{T+1}\tp\f_{T+1}  , & \mbox{if } \hat{\w}\in\bbW^\calC \\
2\mu_{T+1}^2+3\f_{T+1}\tp\f_{T+1}, & \mbox{if } \hat{\w}\in \bbW^\calD \text{ or } \bbW^\calE
\end{cases}.
\end{align*}

\subsection{Unconditional mean squared forecasting error}
\label{sec:UMSFE}
In practice, the candidate forecasts are obtained with errors and thus random, rendering the combination weights also random. This subsection examines the MSFE explicitly accounting for the randomness of the weights and candidate forecasts. In this case, we redefine $\mu_{T+1}=\Exp(y_{T+1})$. The unconditional MSFE of the combined forecast can be written as
\begin{align}
&\quad\Exp(y_{T+1}-\hat{y}_{T+1})^2\notag\\
&=\Exp\left\{y_{T+1}-\mu_{T+1}+\mu_{T+1}-\Exp(\f_{T+1}\tp\hat{\w})+\Exp(\f_{T+1}\tp\hat{\w})-\f_{T+1}\tp\hat{\w}\right\}^2\notag\\
&=\Exp(y_{T+1}-\mu_{T+1})^2+\Exp\left\{\mu_{T+1}-\Exp(\f_{T+1}\tp\hat{\w})\right\}^2+
\Exp\left\{\Exp(\f_{T+1}\tp\hat{\w})-\f_{T+1}\tp\hat{\w}\right\}^2-2\Cov(y_{T+1},\f_{T+1}\tp\hat{\w})\notag\\
&\quad+ 2\Cov\left\{y_{T+1}-\mu_{T+1}, \mu_{T+1}-\Exp(\f_{T+1}\tp\hat{\w})\right\}+2\Cov\left\{\mu_{T+1}-\Exp(\f_{T+1}\tp\hat{\w}), \Exp(\f_{T+1}\tp\hat{\w})-\f_{T+1}\tp\hat{\w}\right\}\notag\\
&=\sigma^2+\left\{\mu_{T+1}-\Exp(\f_{T+1}\tp\hat{\w})\right\}^2+
\Var(\f_{T+1}\tp\hat{\w}) -2\Cov(y_{T+1},\f_{T+1}\tp\hat{\w}),\label{eq:EyT+1}
\end{align}
where the last equality is due to the fact that $\mu_{T+1}-\Exp(\f_{T+1}\tp\hat{\w})$ is nonrandom. The first three terms in~\eqref{eq:EyT+1} are the same as in \eqref{eq:mse} except that all moments are unconditional, while the final and additional covariance term is precisely due to the randomness of $\hat{\w}$ and the fact that both $\hat{\w}$ and $y_{T+1}$ depend on $\f_{T}$. We examine the three terms in turn. The first term $\sigma^2$ is the variance of disturbance that is common across forecasting methods. To calculate the second term, we note that:
\begin{align}\label{eq:Efw}
\Exp(\f_{T+1}\tp\hat{\w})
&=\Exp(\f_{T+1})\tp\Exp(\hat{\w})+\tr \left\{\Cov(\f_{T+1},\hat{\w})\right\}.
\end{align}
Denote $\eta_{\mu_{T+1}}=\mu_{T+1}-\Exp(\f_{T+1})\tp\Exp(\hat{\w})$. Then, by~\eqref{eq:Efw} and Cauchy-Schwarz inequality, we have:
\begin{align}\label{eq:EmuT+12}
&\quad\{\mu_{T+1}-\Exp({\f}_{T+1}\tp\hat{\w})\}^2\notag\\
&=\eta_{\mu_{T+1}}^2-2\eta_{\mu_{T+1}}\tr \left\{\Cov(\f_{T+1},\hat{\w})\right\}+\tr^2 \left\{\Cov(\f_{T+1},\hat{\w})\right\}\notag\\
&\leq 2\eta^2_{\mu_{T+1}}+2\tr^2 \left\{\Cov(\f_{T+1},\hat{\w})\right\}\notag\\
&\leq 2\eta^2_{\mu_{T+1}}+2 \tr\left\{\Var(\f_{T+1})\right\}\tr\left\{\Var(\hat{\w})\right\}\notag\\
&\leq 2\eta^2_{\mu_{T+1}}+2 \tr\left\{\Var(\f_{T+1})\right\}
\Exp(\|\hat{\w}\|^2)\notag\\
&\leq 2\eta^2_{\mu_{T+1}}+2 \tr\left\{\Var(\f_{T+1})\right\} \cdot\begin{cases}
\infty, & \mbox{if } \w\in\bbW^\calA \text{ or } \bbW^\calB\\
S, & \mbox{if } \w\in\bbW^\calC \\
1, & \mbox{if } \w\in\bbW^\calD \text{ or }\bbW^\calE
\end{cases}.
\end{align}
The above inequality suggests that the bias of combined forecasts is bounded by a non-random bias $\eta_{\mu_{T+1}}$ and the variance $\Var(\f_{T+1})$ that depends on the constraints. The constraints determine the upper bound of variance because they affect the variation of random $\hat{\w}$, which further influence the bounds of $\|\Exp(\hat{\w}|{\f}_{T+1})\|^2$ and $\tr\{\Var(\hat{\w}|{\f}_{T+1})\}$.
For the third term in~\eqref{eq:EyT+1}, we can show that:
\begin{align}\label{eq:VarfT+1w}
\quad\Var(\f_{T+1}\tp\hat{\w})
&\leq \Exp(\hat{\w}\tp\f_{T+1}\f_{T+1}\tp\hat{\w})\notag\\
&=\Exp\left\{\Exp(\hat{\w}\tp\f_{T+1}\f_{T+1}\tp\hat{\w}|\f_{T+1})\right\}\notag\\\
&\leq \Exp\left[\Exp(\hat{\w}|\f_{T+1})\tp \f_{T+1}\f_{T+1}\tp\Exp(\hat{\w}|\f_{T+1})+\tr \left\{ \f_{T+1}\f_{T+1}\tp \Var(\hat{\w}|\f_{T+1})\right\}\right]\notag\\
&\leq \Exp\left[\lambda_{\max}(\f_{T+1}\f_{T+1}\tp)\left\|\Exp(\hat{\w}|\f_{T+1})\right\|^2+ \tr(\f_{T+1}\f_{T+1}\tp)\tr \left\{  \Var(\hat{\w}|\f_{T+1})\right\}\right]\notag\\
&= \Exp\left[\|\f_{T+1}\|^2\left\|\Exp(\hat{\w}|\f_{T+1})\right\|^2+ \|\f_{T+1}\|^2\tr \left\{  \Var(\hat{\w}|\f_{T+1})\right\}\right].
\end{align}
This suggests that $\f_{T+1}$ and the first and second-order moments of $\hat{\w}$ play a vital role in the combination variance. Since the distribution of $\f_{T+1}$ is unknown, we cannot analytically derive the variance. Nevertheless, we can examine how the upper bound of \eqref{eq:VarfT+1w} is related to different weight constraints. We note that:
\begin{gather}
\|\Exp(\hat{\w}|{\f}_{T+1})\|^2\leq \begin{cases}
\infty, & \mbox{if } \w\in\bbW^\calA \text{ or } \bbW^\calB\\
S, & \mbox{if } \w\in\bbW^\calC \\
1, & \mbox{if } \w\in\bbW^\calD \text{ or }\bbW^\calE
\end{cases},\\
\intertext{and}
\tr\{\Var(\hat{\w}|{\f}_{T+1})\}\leq \tr\left\{\Exp(\hat{\w}\hat{\w}\tp|{\f}_{T+1})\right\}
\leq \begin{cases}
\infty, & \mbox{if } \w\in\bbW^\calA \text{ or } \bbW^\calB\\
S, & \mbox{if } \w\in\bbW^\calC \\
1, & \mbox{if } \w\in\bbW^\calD \text{ or }\bbW^\calE
\end{cases}.\label{eq:trVarw}
\end{gather}
Hence, combining \eqref{eq:EyT+1} with  \eqref{eq:EmuT+12}--\eqref{eq:trVarw}, we can obtain the upper bound of the MSFE of $\hat{y}_{T+1}$ under various weight constraints as
\begin{align}\label{eq:EyT+1-hyT+11}
&\quad\Exp(y_{T+1}-\hat{y}_{T+1})^2\notag\\
&\leq
\begin{cases}
\infty, & \mbox{if } \w\in\bbW^\calA\text{ or }\bbW^\calB \\
\sigma^2+2\eta^2_{\mu_{T+1}}+2S\tr\left\{\Var(\f_{T+1})\right\}+2S\Exp(\|\f_{T+1}\|^2) -2\Cov(y_{T+1},\f_{T+1}\tp\hat{\w}) , & \mbox{if } \w\in\bbW^\calC \\
\sigma^2+2\eta^2_{\mu_{T+1}}+2\tr\left\{\Var(\f_{T+1})\right\}+2\Exp(\|\f_{T+1}\|^2) -2\Cov(y_{T+1},\f_{T+1}\tp\hat{\w}), & \mbox{if } \w\in\bbW^\calD \text{ or }\bbW^\calE
\end{cases}.
\end{align}
Despite an unknown covariance still appearing in the upper bound, the above bounds suggest that more restrictive constraints, which limit the variation of $\hat{\w}$, reduce the bound of MSFE of the combined forecast.

If the candidate forecasts are unbiased, namely $\Exp(\f_{T+1})=\mu_{T+1}\one$, we have $\eta_{\mu_{T+1}}=\mu_{T+1}-\mu_{T+1}\Exp(\one\tp\hat{\w})=0$, then we have:
\begin{align}\label{eq:EmuT+122}
\quad\{\mu_{T+1}-\Exp({\f}_{T+1}\tp\hat{\w})\}^2
&=\tr^2 \left\{\Cov(\f_{T+1},\hat{\w})\right\}\notag\\
&\leq \tr\left\{\Var(\f_{T+1})\right\}\tr\left\{\Var(\hat{\w})\right\}\notag\\
&\leq  \tr\left\{\Var(\f_{T+1})\right\} \cdot\begin{cases}
\infty, & \mbox{if } \w\in\bbW^\calA \text{ or } \bbW^\calB\\
S, & \mbox{if } \w\in\bbW^\calC \\
1, & \mbox{if } \w\in\bbW^\calD \text{ or }\bbW^\calE
\end{cases}.
\end{align}
Comparing with~\eqref{eq:EmuT+12}, the above inequality shows that when candidate forecasts are unbiased, the bias of the combined forecast has a smaller upper-bound under weight constraints $\bbW^\calC$, $\bbW^\calD$ and $\bbW^\calE$, which further leads to smaller upper bounds of MSFE than \eqref{eq:EyT+1-hyT+11}, that is,
\begin{align}\label{eq:EyT+1-hyT+12}
&\quad\Exp(y_{T+1}-\hat{y}_{T+1})^2\notag\\
&\leq
\begin{cases}
\infty, & \mbox{if } \w\in\bbW^\calA\text{ or }\bbW^\calB \\
\sigma^2+S\tr\left\{\Var(\f_{T+1})\right\}+2S\Exp(\|\f_{T+1}\|^2) -2\Cov(y_{T+1},\f_{T+1}\tp\hat{\w}) , & \mbox{if } \w\in\bbW^\calC \\
\sigma^2+\tr\left\{\Var(\f_{T+1})\right\}+2\Exp(\|\f_{T+1}\|^2) -2\Cov(y_{T+1},\f_{T+1}\tp\hat{\w}), & \mbox{if } \w\in\bbW^\calD \text{ or }\bbW^\calE
\end{cases}.
\end{align}

Furthermore, if the candidate forecasts are all unbiased and uncorrelated with the weights\footnote{This happens, for example, when the forecasts and weights are obtained from different samples.}, namely $\Exp(\f_{T+1}\tp\hat{\w})=\mu_{T+1}$, we have $\{\mu_{T+1}-\Exp({\f}_{T+1}\tp\hat{\w})\}^2=0$, then the combined forecast produces even smaller MSFE than~\eqref{eq:EyT+1-hyT+12} as
\begin{align*}
&\quad\Exp(y_{T+1}-\hat{y}_{T+1})^2\notag\\
&=\sigma^2+
\Var(\f_{T+1}\tp\hat{\w}) -2\Cov(y_{T+1},\f_{T+1}\tp\hat{\w})\notag\\
&\leq
\begin{cases}
\infty, & \mbox{if } \w\in\bbW^\calA\text{ or }\bbW^\calB \\
\sigma^2+2S\Exp(\|\f_{T+1}\|^2) -2\Cov(y_{T+1},\f_{T+1}\tp\hat{\w}) , & \mbox{if } \w\in\bbW^\calC \\
\sigma^2+2\Exp(\|\f_{T+1}\|^2) -2\Cov(y_{T+1},\f_{T+1}\tp\hat{\w}), & \mbox{if } \w\in\bbW^\calD \text{ or }\bbW^\calE
\end{cases}.
\end{align*}
Last, we can observe that the MSFE has a close relationship with the prediction interval. Considering a symmetric prediction interval $[\hat{y}_{T+1}-l, \hat{y}_{T+1}+l]$ for $y_{T+1}$, we have:
\begin{align}
\Pr\left(y_{T+1}\in [\hat{y}_{T+1}-l, \hat{y}_{T+1}+l]\right)
&=\Pr\left(|y_{T+1}-\hat{y}_{T+1}|\leq l\right)\notag\\
&\geq 1- l^{-2} \Exp(y_{T+1}-\hat{y}_{T+1})^2.
\end{align}
Thus, if $l>\sqrt{\alpha^{-1} \Exp(y_{T+1}-\hat{y}_{T+1})^2}$, the coverage probability of $[\hat{y}_{T+1}-l, \hat{y}_{T+1}+l]$ exceeds $1-\alpha$. Furthermore, the minimum length, $l_{\min}$, of the prediction interval tends to be smaller when the $\Exp(y_{T+1}-\hat{y}_{T+1})^2$ is smaller, which can serve as a criterion for determining which weight space is better. An algorithm based on this idea, designed to select weight constraints, is presented in Section \ref{sec:numerical_method}.

\subsection{Uniqueness}
In this subsection, we examine the uniqueness of weights resulting from different weight constraints. Uniqueness is a fundamental property of an optimization problem. The \emph{ex ante} knowledge of uniqueness is helpful to guide us to search for the optimal weights and study the convergence of the weights.

According to convex optimization theory \citep{Stephen2004Convex}, we known that in general the optimal solution of a concave function on a convex set is unique. Note that the weight spaces of $\bbW^\calA, \ldots, \bbW^\calD$ are convex, and an objective function is concave if all eigenvalues of the Hessian matrix are positive. Thus, we shall verify the objective function of each method in order.

First, the regression-based method computes the weights by minimizing the squared loss function, namely $\|\y-\F\w\|^2$. Clearly, when $\lambda_{\min}(T^{-1}\F\tp\F)>0$, the objective function is concave, where $\lambda_{\min}(\cdot)$ represents the smallest eigenvalue. Hence, $\hat{\w}^\calA_\reg, \ldots, \hat{\w}^\calD_\reg$ are all unique if
$\lambda_{\min}(T^{-1}\F\tp\F)>0$. Similarly, $\hat{\w}^{\calA'}_\reg$ is unique when $\lambda_{\min}(T^{-1}\tilde{\F}\tp\tilde{\F})>0$ with $\tilde{\F}=(\one, \F)$.

For optimal averaging-based methods, the objective functions of generalized Mallows and CV defined in~\eqref{eq:objectiveG} and~\eqref{eq:formH} are both of a quadratic form \citep{hansen:2007,hansen.racine:2012}. Thus, $\hat{\w}^\calA_\ma, \ldots, \hat{\w}^\calD_\ma$ are unique when
$
\lambda_{\min}\left(T^{-1}\partial^2 \calD(\w)/\partial \w\tp\partial \w\right)=\lambda_{\min}(T^{-1}\F\tp\F)>0,
$
while $\hat{\w}^\calA_\cv, \ldots, \hat{\w}^\calD_\cv$ are unique when
$
\lambda_{\min}\left(T^{-1}\partial^2 \CV(\w)/\partial \w\tp\partial\w\right)=\lambda_{\min}(T^{-1}\bar{\F}\tp\bar{\F})>0$. 
	
For the weight space $\bbW^\calE$, it is not a convex set but permits a closed-form solution for regression-based and optimal averaging methods. Through the Lagrangian multiplier method, we can obtain the optimal weight $\hat{\w}^\calE_\reg=(\F\tp\F-\nu)^{-1}\F\tp\y$ for the regression-based method, where $\nu$ satisfies $\|(\F\tp\F-\nu)^{-1}\F\tp\y\|=1$. A sufficient condition to guarantee a unique $\hat{\w}^\calE_\reg$ is that $\|(\F\tp\F-\nu)^{-1}\F\tp\y\|>1$ for $\nu\in (\lambda_{\min}(\F\tp\F), \lambda_{\max}(\F\tp\F))$; see Appendix~\hyperref[sec:appendix_B]{B} for the proof.
Similarly, the Mallows and CV averaging weight under $\bbW^\calE$ can be written, respectively, as $\hat{\w}^\calE_{\ma}=-(\F\tp\F-\hat{\nu})^{-1}\bpsi$ and $\hat{\w}^\calE_\cv=(\F\tp\F-\nu)^{-1}\bar{\F}\tp\y$, where $\nu$ satisfies $\|(\F\tp\F-\nu)^{-1}\bpsi\|=1$ in Mallows and $\|(\bar{\F}\tp\bar{\F}-{\nu})^{-1}\bar{\F}\tp\y\|=1$ in CV averaging.
Hence, $\hat{\w}_\ma^\calE$ is unique if $\|(\F\tp\F-\nu)^{-1}\bpsi\|>1$ for $\nu\in (\lambda_{\min}(\F\tp\F),\lambda_{\max}(\F\tp\F))$,  while $\hat{\w}^\calE_\cv$ is unique if $\|(\bar{\F}\tp\bar{\F}-\nu)^{-1}\bar{\F}\tp\y\|>1$ for $\nu\in (\lambda_{\min}(\bar{\F}\tp\bar{\F}),\lambda_{\max}(\bar{\F}\tp\bar{\F}))$.

	By construction, the individual performance-based weights are unique, because they are computed based on a specific performance measure with a one-to-one mapping. Finally, for the eigenvector method with the constraint $\bbW^\calE$, the resulting weight $\hat{\w}^\calE_\eig$ is the eigenvector associated with the smallest eigenvalue of $\M=T^{-1}(\y\otimes \one\tp-\F)\tp(\y\otimes \one\tp-\F)$. Hence, $\hat{\w}^\calE_\eig$ is unique if the smallest eigenvalue of the characteristic polynomial of $(\y\otimes \one\tp-\F)\tp(\y\otimes \one\tp-\F)$ has the multiplicity of one and the first element of $\hat{\w}^\calE_\eig$ is positive.

We summarize the conditions for uniqueness under different weight constraints and estimation methods as follows.

\begin{proposition}\
\begin{itemize}
\item[(1)]If $\lambda_{\min}(T^{-1}\F\tp\F)>0$, $\hat{\w}^\calX_Z$ is unique, where $\calX=\calA$, $\calB$, $\calC$, $\calD$ and $Z$ represents $\reg$ or $\ma$;
\item[(2)]If $\lambda_{\min}(T^{-1}\tilde{\F}\tp\tilde{\F})>0$, $\hat{\w}^{\calA'}_\reg$ is unique, where $\tilde{\F}=(\one,\F)$;
\item[(3)]If $\lambda_{\min}(T^{-1}\bar{\F}\tp\bar{\F})>0$, $\hat{\w}_\cv^\calX$ is unique, where $\calX$ represents $\calA$, $\calB$, $\calC$, $\calD$ and $\bar{\F}=(\f_1^{[-1]\top},\ldots,\f_T^{[-T]\top})\tp$;
\item[(4)]
If $\|(\F\tp\F-\nu)^{-1}\F\tp\y\|>1$ for $\nu\in (\lambda_{\min}(\F\tp\F), \lambda_{\max}(\F\tp\F))$,
then $\hat{\w}^\calE_\reg$ is unique;
\item[(5)]
If $\|(\F\tp\F-\nu)^{-1}\bpsi\|>1$ for $\nu\in (\lambda_{\min}(\F\tp\F),\lambda_{\max}(\F\tp\F))$,
then $\hat{\w}_\ma^\calE$ is unique;
\item[(6)]
If $\|(\bar{\F}\tp\bar{\F}-\nu)^{-1}\bar{\F}\tp\y\|>1$ for $\nu\in (\lambda_{\min}(\bar{\F}\tp\bar{\F}),\lambda_{\max}(\bar{\F}\tp\bar{\F}))$,
then $\hat{\w}^\calE_\cv$ is unique.
\item[(7)]If the smallest eigenvalue of the characteristic polynomial of $(\y\otimes \one\tp-\F)\tp(\y\otimes \one\tp-\F)$ has the multiplicity of 1 and the first element of $\hat{\w}^\calE_\eig$ is positive, then $\hat{\w}^\calE_\eig$ is unique.
\end{itemize}
\end{proposition}

\subsection{Sparsity}
\label{sec:sparsity}
The sparsity of weights is essentially how many elements in the weight vector are zeros, and it is an important target when choosing a set of certain weight constraints. When the number of candidate forecasts is large, and researchers hope to narrow down the candidate forecasts for further examination or interpretation, they typically prefer a sparse weight vector because it suggests that only a few candidate forecasts contribute to the combination. Nevertheless, if the intention is to diversify and take into account as many candidate forecasts as possible for combination, then a dense solution seems a better target. Hence, we analyze the sparsity of weights implied by different weight constraints.

\begin{figure}[H]
\begin{center}
\includegraphics[width=0.5\textwidth]{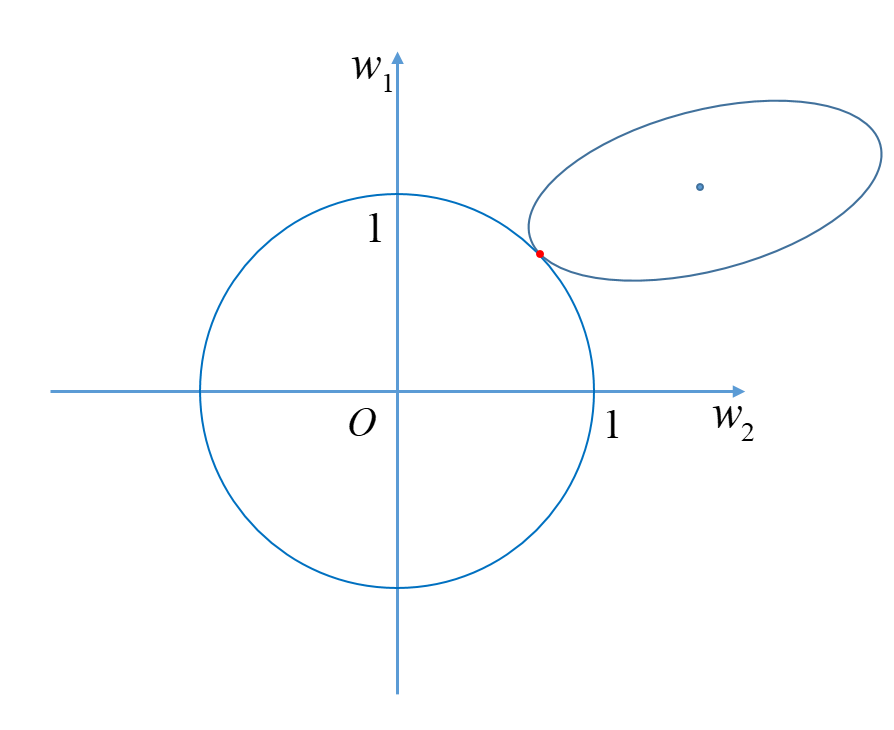}
\end{center}
\caption{The schematic diagram for $\bbW^\calE$. }
\label{fig:eigenvector}
\footnotesize{\emph{Notes:} The circle area represents the space $\bbW^\calE$, the ellipse in the first quadrant represents the equipotential lines of the objective function $f(\x)$.}
\end{figure}

Noting that a quadratic function $f(\x)=\x\tp\A\x+\b\tp\x+c$ is an ellipsoid if and only if $\A$ is a positive definite matrix, and the coordinate of the centre point is $-2^{-1}\A^{-1}\b$. If the solution to the quadratic function lies on the boundaries of coordinate axes, then the solution is sparse.
From the geometrics of weight constraints (see Figure~\ref{fig:regions} for a 2-dimensional example), we can see that under $\bbW^\calA$ the probability of a solution lying on the coordinate axes is zero, so the weights from $\bbW^\calA$, such as $\hat{\w}^\calA_\reg$, $\hat{\w}^\calA_\ma$ and $\hat{\w}^\calA_\cv$, are not sparse.
Similarly, $\bbW^\calB$ and $\bbW^\calE$ do not share boundaries on the coordinate axes, and the feasible solutions lie on the line $\one\tp\w=1$ for $\bbW^\calB$ and $\w\tp\w=1$ for $\hat{\w}^\calE_\reg$ (see Figure~\ref{fig:eigenvector}). Therefore, the resulting weights $\hat{\w}^\calB_\reg$, $\hat{\w}^\calB_\ma$, $\hat{\w}^\calB_\cv$, $\hat{\w}_\eig^\calE$, $\hat{\w}^\calE_\reg$, $\hat{\w}^\calE_\ma$ and $\hat{\w}^\calE_\cv$ are not all sparse.

\begin{figure}[H]
\begin{center}
	\includegraphics[width=0.5\textwidth]{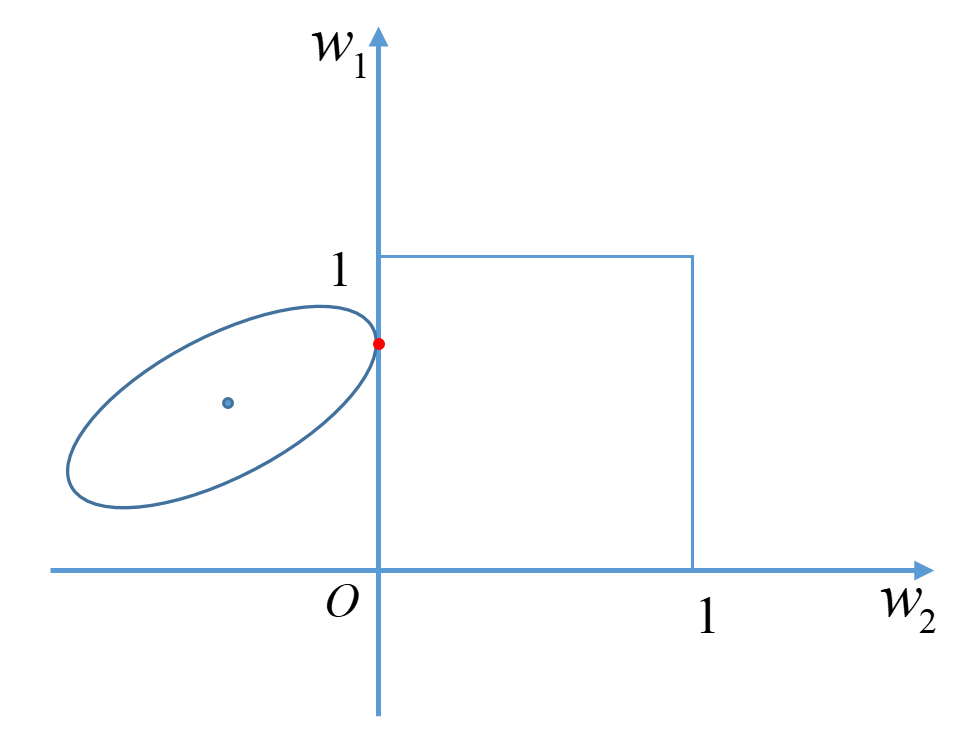}
\end{center}
\caption{The schematic diagram for $\hat{\w}^\calC_\reg$.}\label{fig:solution_D}
\footnotesize{\emph{Notes:} The square area in the first quadrant represents the space $\bbW^\calC$, the ellipse in the second quadrant represents the equipotential lines of the objective function $f(\x)$, and the red interaction point of the two areas on the $y$-axis suggests that $w_2$ is zero.}
\end{figure}

In contrast, the weight space $\bbW^\calC$ and $\bbW^\calD$ contain the boundaries of coordinate axes, and thus sparsity can be achieved under certain conditions. We first investigate the weight space $\bbW^C$. Note that the regression-based method without constraints produces the center point at $-2^{-1}(\F\tp\F)^{-1}\F\tp\y$.
If $\lambda_{\min}(T^{-1}\F\tp\F)>0$, $-2^{-1}(\F\tp\F)^{-1}\F\tp\y\notin \bbW^\calC$ and at least one element of the center points is negative, the resulting weight $\hat{\w}^\calC_\reg$ would be sparse, namely the solution of least-squares optimization reaches the boundaries of coordinate axes, so that some entries of the solution are zeros.
We illustrate this in a 2-dimensional situation in Figure~\ref{fig:solution_D}, where the square area in the first quadrant represents the space $\bbW^\calC$, the ellipse in the second quadrant represents the equipotential lines of the objective function $f(\x)$, and the red interaction point of the two areas on the $y$-axis suggests that $w_2$ is zero. Similarly, the Mallows' averaging weight $\hat{\w}^\calC_\ma$ is sparse if $\lambda_{\min}(T^{-1}\F\tp\F)>0$ and $-(\F\tp\F)^{-1}\F\tp\bpsi\notin \bbW^\calC$ with at least one element being negative. The CV averaging weight $\hat{\w}^\calC_\cv$ is sparse if $\lambda_{\min}(T^{-1}\bar{\F}\tp\bar{\F})>0$ and $-(\bar{\F}\tp\bar{\F})^{-1}\bar{\F}\tp\y\notin \bbW^\calC$ with at least one element being negative. The feature of sparse weights for model averaging methods is also discussed by \cite{Yang2020On}.

\begin{figure}[H]
	\begin{center}
		\begin{subfigure}[pic1]{0.4\textwidth}
			\includegraphics[width=\textwidth]{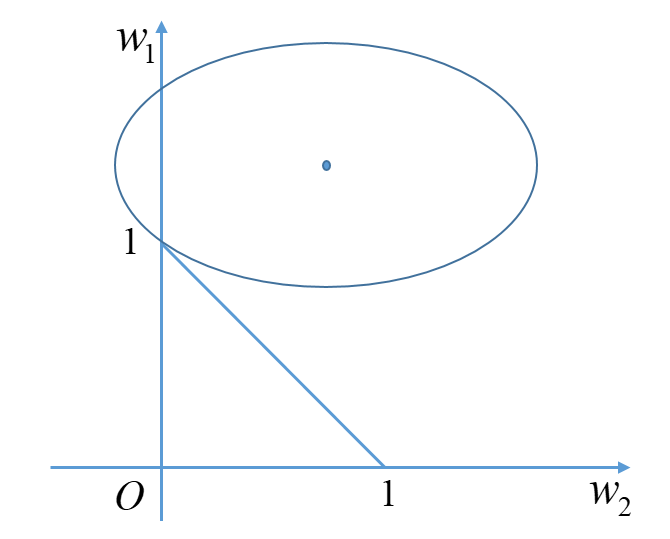}
			\caption{$\{\w|\w\in [0,1]^2, w_1+w_2=1\}$}
		\end{subfigure}
		\begin{subfigure}[pic2]{0.4\textwidth}
			\includegraphics[width=0.85\textwidth]{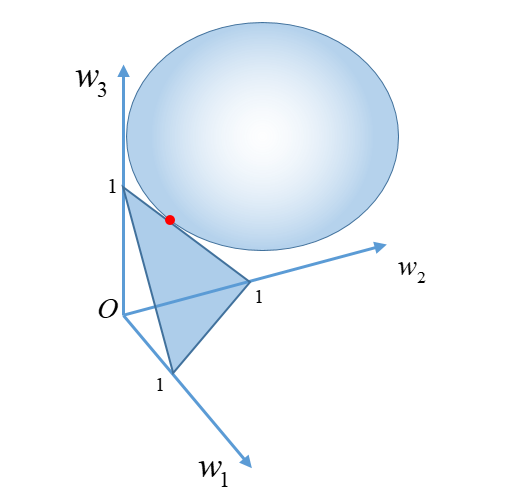}
			\caption{$\{\w|\w\in [0,1]^3, w_1+w_2+w_3=1\}$}
		\end{subfigure}
	\end{center}
	\caption{The schematic diagram for $\hat{\w}^\calD_\reg$.}
	\label{fig:weights_calcualting}
	\footnotesize{\emph{Notes:} The line segment from (0,1) to (1,0) in Figure \ref{fig:weights_calcualting}(a) and the shadow triangle in Figure \ref{fig:weights_calcualting}(b) represent the feasible region of $\w$. The ellipse in Figure~\ref{fig:weights_calcualting}(a) and the ellipsoid in Figure~\ref{fig:weights_calcualting}(b) are the contour line/surface of the objective function. The solution under $\hat{\w}^\calD_\reg$ is the intersection point between the feasible region $\bbW^\calD$ and the contour line/surface line with the smallest distance.}
\end{figure}

Next, we examine the weight space $\bbW^\calD$.
For regression-based methods, a sufficient condition for $\hat{\w}^\calD_\reg$ to be sparse (with probability one) is that it is a boundary point of $[0,1]^S$ but not the tangent point of the plane $\one\tp\w=1$; in other words, from the Kuhn–Tucker condition, the sufficient condition implies that there is not a nonzero constant $\rho_0$ satisfying
\begin{align}\label{eq:61}
-\frac{2}{T}\sum_{t=1}^{T}(y_t-\f_t\tp\hat{\w}^\calD_\reg)\f_t=\rho_0\one.
\end{align}
This condition is illustrated in Figures~\ref{fig:weights_calcualting}(a) and \ref{fig:weights_calcualting}(b) for a 2- and 3-dimensional case, respectively.
For model averaging methods, we can follow similar reasoning to conclude that the weights $\hat{\w}^\calD_\ma$ and $\hat{\w}^\calD_\cv$ are sparse, if there do not exist any nonzero $\rho_1$ and $\rho_2$ satisfying
$\F\tp\F\w+\bpsi=\rho_1\one$ for Mallows and $\bar{\F}\tp\bar{\F}\w-\bar{\F}\tp\y=\rho_2\one$ for CV averaging.

	Finally, the individual performance-based weights lie in the space of $\calD$ by construction. Typically, they are not sparse with probability being 1, because the performance measure of a candidate forecast is often nonzero, and the weights are also normalized.

We summarize the conditions of sparsity for different weight constraints and estimation methods in the following proposition and Table~\ref{tab:sparse}.
\begin{proposition}\
\begin{itemize}
	\item[(1)] The optimal weights $\hat{\w}^\calA_\reg$, $\hat{\w}^{\calA'}_\reg$, $\hat{\w}^\calB_\reg$, $\hat{\w}^\calE_\reg$,
	$\hat{\w}^\calA_\ma$, $\hat{\w}^\calB_\ma$, $\hat{\w}^\calE_\ma$, $\hat{\w}^\calA_\cv$, $\hat{\w}^\calB_\cv$, $\hat{\w}^\calE_\cv$, $\hat{\w}^\calD_\pf$ and $\hat{\w}^\calE_\eig$ usually are not sparse.
	\item[(2)] If $\lambda_{\min}(T^{-1}\F\tp\F)>0$, $-2^{-1}(\F\tp\F)^{-1}\F\tp\y\notin \bbW^\calC$ and there exists a vector $\e_i$ for $i=1, 2, \ldots, S$ such that $\e_i\tp(\F\tp\F)^{-1}\F\tp\y>0$, then optimal weight $\hat{\w}^\calC_\reg$ is sparse, where the $i$th entry of $\e_i$ is 1 and others are 0.
	\item[(3)] If $\lambda_{\min}(T^{-1}\F\tp\F)>0$, $-(\F\tp\F)^{-1}\F\tp\bpsi\notin \bbW^\calC$ and there exists a vector $\e_i$ for $i=1, 2, \ldots, S$ such that $\e_i\tp(\F\tp\F)^{-1}\F\tp\bpsi>0$, the optimal weight $\hat{\w}^\calC_\ma$ is sparse, where the $i$th entry of $\e_i$ is 1 and others are 0.
	\item[(4)]If $\lambda_{\min}(T^{-1}\bar{\F}\tp\bar{\F})>0$, $-(\bar{\F}\tp\bar{\F})^{-1}\bar{\F}\tp\y\notin \bbW^\calC$ and there exists a vector $\e_i$ for $i=1, 2, \ldots, S$ such that $\e_i\tp(\bar{\F}\tp\bar{\F})^{-1}\bar{\F}\tp\y>0$, the optimal weight $\hat{\w}^\calC_\cv$ is sparse, where the $i$th entry of $\e_i$ is 1 and others are 0.
	\item[(5)]If $\lambda_{\min}(T^{-1}\F\tp\F)>0$ and there is not a nonzero constant $\rho_0$ satisfying $T^{-1}\sum_{t=1}^{T}(y_t-\f_t\tp\hat{\w}^\calD_\pf)\f_t=\rho_0\one$, then the solution of weight $\hat{\w}^\calD_\reg$ is sparse.
	\item[(6)] If $\lambda_{\min}(T^{-1}\F\tp\F)>0$ and there is not a nonzero constant $\rho_1$ satisfying $\F\tp\F\w+\bpsi=\rho_1 \one$, then the solution of weight $\hat{\w}^\calD_\ma$ is sparse.
\item[(7)] If $\lambda_{\min}(T^{-1}\bar{\F}\tp\bar{\F})>0$ and there is not a nonzero constant $\rho_2$ satisfying $\bar{\F}\tp\bar{\F}\w-\bar{\F}\tp\y=\rho_2\one$, then the solution of weight $\hat{\w}^\calD_\cv$ is sparse.
\end{itemize}
\end{proposition}

\begin{table}[H]
\caption{The property for sparseness.}\label{tab:sparse}
\centering
\begin{tabular}{c|c|p{1.4cm}<{\centering}|p{1.1cm}<{\centering}|c|c}
	\toprule
	\multirow{2}{*}{regions} &
	regression-based & \multicolumn{2}{c|}{model averaging} & individual performance-based & eigenvector\\
	\cline{2-6}
	&   $\reg$ &  $\ma$  & $\cv$   &   $\pf$ & $\eig$\\
	\midrule
	$\calA$ &  $\times$ & $\times$ & $\times$  & --- & --- \\
	$\calB$ &  $\times$ & $\times$ & $\times$ & --- & ---  \\
	$\calC$ &  $\surd$ & $\surd$ & $\surd$     &  ---    &   ---    \\
	$\calD$ &  $\surd$ & $\surd$ & $\surd$ &   $\times$  &     ---  \\
	$\calE$ &   $\times$& $\times$ & $\times$ &  --- & $\times$   \\
	\bottomrule
\end{tabular}\\
{
	\begin{flushleft}
		\small Note: ``$\surd$" indicates that the weight is sparse under some conditions; ``$\times$" indicates that the weight is not sparse with probability equal to 1; ``---" means that the case is ambiguous or does not exist.
	\end{flushleft}
}
\end{table}

\section{The guidance to select a proper weight space}\label{sec:guidance}
\subsection{From the Bayesian perspective}\label{sec:Bayesian}
On one hand,
the weights in $\bbW^\calC$ and $\bbW^\calD$ are recommended because they resemble probabilities.
Bayesian model averaging (BMA) combines forecasts based on the posterior probability assigned to their associated models, and thus the weights of BMA fall into the space of $[0,1]^S$, namely they belong to either $\bbW^\calC$ or $\bbW^\calD$. More specifically,
consider forecasts obtained from two models, labeled as ``$\model_1$" and ``$\model_2$", BMA obtains the forecast from the unconditional mean as
\begin{align}
\Exp(y_{T+1})=\Pr(\model_1){\Exp}(y_{T+1}|\model_1)+\Pr(\model_2){\Exp}(y_{T+1}|\model_2),
\end{align}
where $\Pr(\model_i)$ denotes the probability that $\model_i$ coincides with the data generating process and ${\Exp}(y_{T+1}|\model_i)$ is the conditional expectation of $y_{T+1}$ given $\model_i$ for $i=1,2$.
Thus, the posterior probability $\Pr(\model_1$) and $\Pr(\model_2$) serve as weights in BMA.

On the other hand, the weight spaces are supported by their corresponding prior distributions.
Consider the weights as random variables with a density given by $p(\w)=g(\w)1_{\w\in \bbW}(\w)$. The distribution of $y_t$ conditional on $\f_t$ and $\w$ is $p(y_t|\f_t,\w)=\Nor(\f_t\tp\w, 1|\f_t,\w)\propto \exp\left\{-(y_t-\f_t\tp\w)^2 /2\right\}$. Thus, its posterior distribution is
$p(\w|\f_t,y_t)\propto p(y_t|\f_t,\w)p(\w)\propto \exp[-(y_t-\f_t\tp\w)^2 /2+\log\{g(\w)\}]1_{\w\in\bbW}(\w)$. In this case, the maximum a posteriori (MAP) estimator of $\w$ is
\begin{align}\label{eq:map}
&\argmax_\w \prod_{t=1}^{T}\exp\left[-2^{-1}(y_t-\f_t\tp\w)^2+\log\{g(\w)\}\right] 1_{\w\in\bbW}(\w)\notag\\
&=\argmin_{\w\in\bbW}\frac{1}{T}\sum_{t=1}^{T}(y_t-\f_t\tp\w)^2-\log\{g(\w)\}.
\end{align}
From \eqref{eq:map}, we know the weight space refers to the support of some prior density.

For example, if the prior density $g(\w)$ is an $S$-dimensional normal distribution, then $\bbW^\calA$ is a better choice. If $g(\w)\propto \exp\{-\sums(w_s-0.5)^2\}$ is an $(S-1)$-dimensional normal distribution in $\bbW^B$, then $\bbW^B$ is a better choice. For the bounded regions, we can simply consider the uniform distribution as the prior distribution. For instance, $g(\w)=1_{\w\in\bbW^\calC}(\w)$ for $\bbW^\calC$, $g(\w)=1_{\w\in\bbW^\calD}(\w)/m(\bbW^\calD)$ for $\bbW^\calD$ and $g(\w)=1_{\w\in\bbW^\calE}(\w)/m(\bbW^\calE)$ for $\bbW^\calE$, where $m(\cdot)$ represents the cardinality for a set.\footnote{For example, $m(\bbW^\calD)=\sqrt{2}$ is the length of a segment in two-dimensional space, and $m(\bbW^\calD)=\sqrt{3}/2$ is the area of a triangle in three-dimensional space; $m(\bbW^\calE)=2\pi$ is the perimeter of a circle in two-dimensional space; and $m(\bbW^\calE)=4\pi$ is the area of a sphere in three-dimensional space.} In particular, for $\bbW^\calD$, we can also consider the prior distribution to be an $S$-dimensional Dirichlet distribution:
\begin{gather}
p(\w|\balpha)=\frac{\Gamma(\sum_{s=1}^{S} \alpha_i)}{\prod_{s=1}^{S}\Gamma(\alpha_i)}w_1^{\alpha_1-1}\cdots w_S^{\alpha_S-1}1_{\w\in\bbW^\calD},
\end{gather}
where the parameter $\balpha$ is an $S$-vector with components $\alpha_s>0$, $\Gamma(x)$ is the Gamma function. When $\alpha_s=1$ for $s=1, \ldots, S$, this distribution degenerates to $1_{\w\in\bbW^\calD}(\w)/m(\bbW^\calD)$.

\begin{remark}
Based on \eqref{eq:map}, if we consider the uniform distribution for bounded regions, we find that the weight constraints lead to the penalties on optimization process:
\begin{gather*}
\min_{\w} \|y_t\one-\f_t\tp\w\|^2-\bmu\tp \w-\bnu\tp (\one-\w), \text{ for } \bbW^\calC,\\
\min_{\w} \|y_t\one-\f_t\tp\w\|^2+\lambda \w\tp \one-\bmu\tp \w,\text{ for } \bbW^\calD,\\
\min_{\w} \|y_t\one-\f_t\tp\w\|^2+\lambda \w\tp\w,\text{ for } \bbW^\calE,
\end{gather*}
where $\lambda$, $\bmu$, $\bnu$ are the lagrangian multipliers for exact optimal solution, and $\lambda$, $\bmu$, $\bnu$ are predefined some positive numbers  for soft constraints.
\end{remark}

\subsection{A numerical method to choose weight space}
\label{sec:numerical_method}
From Subsection \ref{sec:UMSFE}, although different weight spaces have different influences on variance and bias, they collectively impact the predictions. Therefore, we aim to use the length of the prediction interval as a criterion for selecting an appropriate  weight space. In this context, we employ the technique of conformal inference\citep{Jing2018Distribution,Yang2025selection} to obtain numerical results of the prediction interval, which will guide our selection of the weight space.  This idea is summarized in Algorithm \ref{alg:select-weightregion}.

\begin{algorithm}[H]
\label{alg:select-weightregion}
\caption{Selecting weight constraints by conformal inference.}
\SetAlgoLined
\SetKwFunction{Model}{Model}
\SetKwInOut{Input}{Input}
\SetKwInOut{Output}{Output}
\SetKw{IFor}{for}
\SetKw{Iin}{in}
\Input{$\{(\x_t, y_t)\}_{t=1}^T$, miscoverage level $\alpha$}
\Output{$\bbW$}
Randomly split $\{1, \ldots, T\}$ into three equal-sized subsets $\calI_1, \calI_2, \calI_3$\;
\For{$\calX$ \Iin $\{\calA, \calB, \ldots, \calE\}$}{
$\hat{f}_{(1)}, \ldots, \hat{f}_{(S)}=\Model(\{(\x_t, y_t): t\in\calI_1\})$\;\tcc{$\Model(\cdot)$ means the module to train $S$ candidate models}
$\hat{\w}=\argmin_{\w\in\bbW^\calX}\sum_{t\in\calI_2} \left\{\sums w_s\hat{f}_\s(\x_t)-y_t\right\}^2$\;
$R_t=|y_t-\sums\hat{w}_s\hat{f}_\s(\x_t)|$ \IFor $t\in\calI_3$ \;
$l^\calX$=the $k$-th smallest value in $\{R_i: i\in\calI_3\}$, where $k=\lceil (n+1)(1-\alpha)\rceil$
}
$\calX=\argmin_\calX l^\calX$\;
$\bbW=\bbW^\calX$\;
\end{algorithm}

\section{Simulation}\label{sec:simulation}
	This section numerically verifies the properties of estimated weights obtained from different constraints and methods via a simulation study. We consider the following data generating process (DGP):
	\begin{align*}
	y_t=\x_t\tp\bbeta+\epsilon_t, \quad t=1, 2, \ldots T,
	\end{align*}
	where $\bbeta$ is a $p$-dimensional vector, $\epsilon_t$ is independently drawn from a standard normal distribution. We consider four cases of regressors with distinct correlations and distributions:
	
	\begin{itemize}[itemsep=0pt,topsep=0pt,partopsep=0pt]
		\item[] \textbf{Case 1}: $\x_t\sim \Nor(\boldsymbol{0}, \bSigma)$, where $\bSigma=\I_{p\times p}$.
		\item[] \textbf{Case 2}: $\x_t\sim \Nor(\boldsymbol{0}, \bSigma)$, where $\bSigma=(0.7^{|i-j|})_{p\times p}$.
		\item[] \textbf{Case 3}: $\x_t$ follows a multivariate $t$ distribution with the location vector $\boldsymbol{0}$, the scale matrix $\bSigma=\I_{p\times p}$ (note that $\bSigma\neq \Cov(\x)$), and the degree of freedom $\nu=2$.
		\item[] \textbf{Case 4}: $\x_t$ follows a multivariate $t$ distribution with the location vector $\boldsymbol{0}$, the scale matrix $\bSigma=(0.7^{|i-j|})_{p\times p}$, and the degree of freedom $\nu=2$.
	\end{itemize}
	Cases 1 and 2 consider normally distributed regressors, while Cases 3 and 4 consider regressors with a flatter tail. The regressors are correlated with each other in Cases 2 and 4 but not in Cases 1 and 3. We shall examine how the distribution and correlation influence the relation among candidate models and further the optimal weights.
	
	To examine how the quality of candidate models affects the weight optimization, we also consider four ways to construct the candidate models, which ultimately differ in the sets of regressors included in the model.
	\begin{itemize}[itemsep=0pt,topsep=0pt,partopsep=0pt]
		\item[] \textbf{Set 1}:
		The covariates of the $s$th candidate model is $\x_{t}^{(s)}=(x_{t,4(s-1)+1}, \ldots, x_{t,\min(4s,d)})\tp$ for $s=1, 2, \ldots, \lceil d/4\rceil$, $t=1, \ldots, T$, where $d$ determines the number of regressors.
		\item[] \textbf{Set 2}: The same set as above except excluding the last two regressors, namely $\x_{t}^{(s)}=(x_{t,4(s-1)+1}, \ldots, x_{t,\min(4s,d-2)})\tp$ for $s=1, 2, \ldots, \lceil (d-2)/4\rceil$, $t=1, \ldots, T$.
		\item[] \textbf{Set 3}: $\x_{t}^{(s)}=(x_{t,s+2}, \ldots, x_{t,\min(s+4,d)})\tp$ for $s=1, 2, \ldots, \lceil d/4\rceil$, $t=1, \ldots, T$.
		\item[] \textbf{Set 4}: $\x_{t(s)}=(x_{t,s+2}, \ldots, x_{t,\min(s+4,d-2)})\tp$ for $s=1, 2, \ldots, \lceil (d-2)/4\rceil$, $t=1, \ldots, T$.
	\end{itemize}
	Note that regressors in Sets 1 and 2 do not overlap, such that the candidate forecasts are less correlated. In contrast, Sets 3 and 4 allow candidate models to share regressors, leading to a higher correlation between candidate forecasts. Sets 2 and 4 intentionally omit some regressors, such that all candidate models are misspecified.
We set $T=10000$ and $d=42$, and consider 16 scenarios (4 Cases $\times$ 4 Sets ).

	Table~\ref{tab:SSR} presents the SSR for different weight constraints and estimation methods. First, we find that $\calA$ generally produces the lowest SSR, while the weights obtained from $\calD$ are associated with the largest SSR, confirming the theory of Section \ref{sec:SSR}, that is, a region with a larger range tends to result in a lower SSR. We also note that the SSR of regression-based and model-averaging methods is comparable and lower than that of other methods. This result is mainly because these two methods both minimize the quadratic loss of residuals or its approximation.
	
	Table~\ref{tab:biasedness} presents the empirical biasedness. Thanks to the inclusion of an intercept, $\hat{\w}^{\calA'}_\reg$ leads to an unbiased combined forecast, confirmed by the first column of the table. Other methods are generally biased except when candidate forecasts are unbiased and the sum-to-unity constraint is imposed.
	
	Next, we evaluate the MSFE using the test sample and present the results in Table \ref{tab:MSFE}. We find that the MSFE of the unconstrained or less constrained combined forecast ($\w\in\bbW^\calA$) is generally smaller than those of (more) constrained combination. However, in Cases 3-4 and Sets 3-4, the MSFE resulting from $\bbW^\calA$ and $\bbW^\calB$ is larger than that from $\bbW^\calD$. This is because in these ``difficult'' cases to forecast, a larger forecasting variance is expected and less restricted weights may also lead to overfitting. On the contrary, more regularization in the constraint reduces the variance and helps avoid overfitting, albeit at the cost of sacrificing some bias.

	Finally, to examine the sparsity property, we report the percentage of zeros in the resulting weight vector under different estimation methods and constraints in Tables \ref{tab:weight_precent_zero}. It shows that $\bbW^\calC$ and $\bbW^\calD$ do result in a large degree of sparsity with many zero elements in the weight vector, as analyzed in Section \ref{sec:sparsity}. In contrast, $\bbW^\calA$, $\bbW^\calB$ and $\bbW^\calE$ do not share boundaries with the coordinates, rendering nonsparse weight vectors.

\begin{sidewaystable}[H]
\caption{Simulation results: Sum of squared residuals ($\times 10^4$).}
\label{tab:SSR}%
\footnotesize
\begin{tabular}{cc|*{13}{c}}
\toprule
Case & Set & $\hat{\w}^\calA_{\reg}$ & $\hat{\w}^\calB_{\reg}$ & $\hat{\w}^\calC_{\reg}$ & $\hat{\w}^\calD_{\reg}$ & $\hat{\w}^\calE_{\reg}$ & $\hat{\w}^\calA_{\ma}$ & $\hat{\w}^\calB_{\ma}$ & $\hat{\w}^\calC_{\ma}$ & $\hat{\w}^\calD_{\ma}$ & $\hat{\w}^\calE_{\ma}$ & $\hat{w}_{\pf\_\SAIC}^\calD $ & $\hat{w}_{\pf\_\SBIC}^\calD$ & $\hat{\w}_\eig^\calE$ \\
\midrule
1     & 1     & 1.006  & 1.036  & 1.008  & 6.097  & 3.803  & 1.006  & 1.036  & 1.009  & 6.097  & 3.803  & 6.857  & 7.956  & 12.614  \\
1     & 2     & 1.490  & 1.514  & 1.492  & 6.097  & 3.920  & 1.490  & 1.514  & 1.492  & 6.097  & 3.920  & 6.857  & 7.956  & 4.613  \\
1     & 3     & 1.490  & 1.499  & 1.493  & 6.504  & 2.816  & 1.492  & 1.499  & 1.494  & 6.504  & 2.816  & 7.093  & 7.956  & 14.938  \\
1     & 4     & 1.978  & 1.987  & 1.979  & 6.509  & 3.117  & 1.982  & 1.987  & 1.979  & 6.509  & 3.117  & 7.093  & 7.956  & 14.509  \\

2     & 1     & 1.117  & 1.214  & 1.118  & 2.096  & 1.650  & 1.117  & 1.214  & 1.118  & 2.096  & 1.650  & 2.124  & 2.438  & 3.419  \\
2     & 2     & 1.225  & 1.304  & 1.227  & 2.096  & 1.691  & 1.226  & 1.304  & 1.227  & 2.096  & 1.691  & 2.124  & 2.438  & 3.366  \\
2     & 3     & 1.123  & 1.215  & 1.163  & 2.148  & 1.560  & 1.123  & 1.215  & 1.163  & 2.148  & 1.560  & 2.252  & 2.439  & 3.641  \\

2     & 4     & 1.261  & 1.335  & 1.292  & 2.173  & 1.630  & 1.261  & 1.335  & 1.293  & 2.173  & 1.630  & 2.202  & 2.439  & 3.565  \\

3     & 1     & 34.834  & 34.919  & 36.165  & 48.752  & 35.232  & 34.834  & 34.919  & 36.165  & 48.752  & 35.232  & 63.946  & 154.732  & 800.374  \\
3     & 2     & 34.849  & 38.335  & 36.165  & 48.752  & 35.768  & 34.849  & 38.335  & 36.165  & 48.752  & 35.768  & 63.946  & 126.125  & 896.235  \\
3     & 3     & 7.967  & 9.013  & 34.415  & 48.025  & 25.766  & 7.967  & 9.013  & 34.415  & 48.025  & 25.766  & 64.828  & 145.688  & 1705.625  \\

3     & 4     & 8.703  & 9.616  & 34.579  & 48.025  & 26.742  & 8.703  & 9.616  & 34.579  & 48.025  & 26.742  & 64.828  & 145.688  & 1638.297  \\

4     & 1     & 8.530  & 8.828  & 9.096  & 11.555  & 8.576  & 8.530  & 8.828  & 9.096  & 11.555  & 8.577  & 15.131  & 27.859  & 159.476  \\
4     & 2     & 8.531  & 9.476  & 9.098  & 11.555  & 8.580  & 8.531  & 9.476  & 9.098  & 11.555  & 8.580  & 15.131  & 23.440  & 168.111  \\
4     & 3    & 2.682  & 2.712  & 9.373  & 12.136  & 7.417  & 2.682  & 2.712  & 9.373  & 12.136  & 7.417  & 15.021  & 27.145  & 349.425  \\

4     & 4     & 3.269  & 3.347  & 9.430  & 12.136  & 7.598  & 3.269  & 3.347  & 9.430  & 12.136  & 7.598  & 15.021  & 27.145  & 333.227  \\

\bottomrule
\end{tabular}
\end{sidewaystable}

\begin{sidewaystable}[H]
\centering
\caption{Simulation results: Empirical biasedness ($\times 10^4$).}
\label{tab:biasedness}
\footnotesize
\begin{tabular}{cc|*{14}{c}}
\toprule
Case & Set & $\hat{\w}^{\calA'}_{\reg}$ & $\hat{\w}^\calA_{\reg}$ & $\hat{\w}^\calB_{\reg}$ & $\hat{\w}^\calC_{\reg}$ & $\hat{\w}^\calD_{\reg}$ & $\hat{\w}^\calE_{\reg}$ & $\hat{\w}^\calA_{\ma}$ & $\hat{\w}^\calB_{\ma}$ & $\hat{\w}^\calC_{\ma}$ & $\hat{\w}^\calD_{\ma}$ & $\hat{\w}^\calE_{\ma}$  & $\hat{w}_{\pf\_\SAIC}^\calD $ & $\hat{w}_{\pf\_\SBIC}^\calD$ & $\hat{\w}_\eig^\calE$ \\
\midrule
1     & 1  &$-3\times 10^{-16}$   & $-0.014 $  &   $-0.021 $  &  $-0.010 $   &  $ 0.000 $ &  $-0.010$  &  $-0.015 $ &  $-0.021$  &  $-0.012$  &  $ 0.000 $ &  $-0.010$  &  $-0.040$  &  $-0.027$  &  $-0.036 $ \\
1     & 2  &$-3\times 10^{-16}$   & $-0.013 $  &   $-0.019 $  &  $-0.010 $   &  $ 0.000 $ &  $-0.009$  &  $-0.014 $ &  $-0.019$  &  $-0.012$  &  $ 0.000 $ &  $-0.009$  &  $-0.040$  &  $-0.027$  &  $-0.016 $ \\
1     & 3  &$ 1\times 10^{-15}$   & $0.046  $  &   $0.034  $  &  $ 0.055 $   &  $ 0.000 $ &  $ 0.051$  &  $ 0.040 $ &  $ 0.034$  &  $ 0.051$  &  $ 0.000 $ &  $ 0.051$  &  $-0.010$  &  $-0.024$  &  $-0.097 $ \\
1     & 4  &$ 3\times 10^{-15}$   & $0.053  $  &   $0.041  $  &  $ 0.055 $   &  $-0.002 $ &  $ 0.054$  &  $ 0.045 $ &  $ 0.041$  &  $ 0.051$  &  $-0.002 $ &  $ 0.054$  &  $-0.010$  &  $-0.024$  &  $-0.098 $ \\
2     & 1  &$-1\times 10^{-15}$   & $0.104  $  &   $0.130  $  &  $ 0.106 $   &  $ 0.036 $ &  $ 0.066$  &  $ 0.104 $ &  $ 0.130$  &  $ 0.106$  &  $ 0.036 $ &  $ 0.066$  &  $ 0.000$  &  $ 0.048$  &  $ 0.046 $ \\
2     & 2  &$-6\times 10^{-16}$   & $0.079  $  &   $0.105  $  &  $ 0.081 $   &  $ 0.036 $ &  $ 0.058$  &  $ 0.080 $ &  $ 0.105$  &  $ 0.081$  &  $ 0.036 $ &  $ 0.058$  &  $ 0.000$  &  $ 0.048$  &  $ 0.055 $ \\
2     & 3  &$ 2\times 10^{-15}$   & $0.102  $  &   $0.114  $  &  $ 0.126 $   &  $ 0.070 $ &  $ 0.094$  &  $ 0.103 $ &  $ 0.114$  &  $ 0.126$  &  $ 0.070 $ &  $ 0.094$  &  $ 0.089$  &  $ 0.054$  &  $ 0.039 $ \\
2     & 4  &$-6\times 10^{-17}$   & $0.095  $  &   $0.100  $  &  $ 0.116 $   &  $ 0.062 $ &  $ 0.084$  &  $ 0.095 $ &  $ 0.100$  &  $ 0.116$  &  $ 0.062 $ &  $ 0.083$  &  $ 0.057$  &  $ 0.054$  &  $ 0.046 $ \\
3     & 1  &$-2\times 10^{-15}$   & $0.322  $  &   $0.341  $  &  $ 0.224 $   &  $ 0.552 $ &  $ 0.291$  &  $ 0.323 $ &  $ 0.341$  &  $ 0.224$  &  $ 0.552 $ &  $ 0.291$  &  $ 0.485$  &  $ 1.428$  &  $ 3.462 $ \\
3     & 2  &$ 4\times 10^{-15}$   & $0.317  $  &   $0.470  $  &  $ 0.224 $   &  $ 0.552 $ &  $ 0.301$  &  $ 0.317 $ &  $ 0.470$  &  $ 0.224$  &  $ 0.552 $ &  $ 0.301$  &  $ 0.485$  &  $ 1.121$  &  $ 3.650 $ \\
3     & 3  &$ 3\times 10^{-12}$   & $0.162  $  &   $0.307  $  &  $ 0.108 $   &  $ 0.476 $ &  $ 0.314$  &  $ 0.161 $ &  $ 0.307$  &  $ 0.108$  &  $ 0.476 $ &  $ 0.314$  &  $ 0.622$  &  $ 1.491$  &  $ 4.745 $ \\
3     & 4  &$-4\times 10^{-13}$   & $0.283  $  &   $0.409  $  &  $ 0.116 $   &  $ 0.476 $ &  $ 0.280$  &  $ 0.282 $ &  $ 0.409$  &  $ 0.116$  &  $ 0.476 $ &  $ 0.280$  &  $ 0.622$  &  $ 1.491$  &  $ 4.609 $ \\
4     & 1  &$-2\times 10^{-15}$   & $0.061  $  &   $0.090  $  &  $ 0.004 $   &  $ 0.217 $ &  $ 0.065$  &  $ 0.061 $ &  $ 0.090$  &  $ 0.004$  &  $ 0.217 $ &  $ 0.065$  &  $ 0.211$  &  $ 0.600$  &  $ 1.755 $ \\
4     & 2  &$-3\times 10^{-15}$   & $0.062  $  &   $0.244  $  &  $ 0.003 $   &  $ 0.217 $ &  $ 0.072$  &  $ 0.062 $ &  $ 0.244$  &  $ 0.003$  &  $ 0.217 $ &  $ 0.072$  &  $ 0.211$  &  $ 0.381$  &  $ 1.826 $ \\
4     & 3  &$ 5\times 10^{-12}$   & $0.122  $  &   $0.107  $  &  $ 0.057 $   &  $ 0.229 $ &  $ 0.173$  &  $ 0.122 $ &  $ 0.107$  &  $ 0.057$  &  $ 0.229 $ &  $ 0.173$  &  $ 0.159$  &  $ 0.312$  &  $ 2.376 $ \\
4     & 4  &$ 3\times 10^{-12}$   & $0.111  $  &   $0.079  $  &  $ 0.058 $   &  $ 0.229 $ &  $ 0.158$  &  $ 0.111 $ &  $ 0.079$  &  $ 0.058$  &  $ 0.229 $ &  $ 0.158$  &  $ 0.159$  &  $ 0.312$  &  $ 2.331 $ \\

\bottomrule
\end{tabular}
\end{sidewaystable}

\begin{sidewaystable}[H]
\caption{Simulation results: Mean squared forecast error.}
\label{tab:MSFE}\footnotesize
\begin{tabular}{cc|*{13}{c}}
\toprule
Case & Set & $\hat{\w}^\calA_{\reg}$ & $\hat{\w}^\calB_{\reg}$ & $\hat{\w}^\calC_{\reg}$ & $\hat{\w}^\calD_{\reg}$ & $\hat{\w}^\calE_{\reg}$ & $\hat{\w}^\calA_{\ma}$ & $\hat{\w}^\calB_{\ma}$ & $\hat{\w}^\calC_{\ma}$ & $\hat{\w}^\calD_{\ma}$ & $\hat{\w}^\calE_{\ma}$ & $\hat{\w}^\calA_{\pf}$ & $\hat{\w}^\calB_{\pf}$ & $\hat{\w}^\calE_{\eig}$ \\
\midrule
1     & 1     & 1.024  & 1.057  & 1.023  & 6.176  & 3.854  & 1.024  & 1.057  & 1.023  & 6.176  & 3.854  & 6.979  & 8.002  & 12.660  \\
1     & 2     & 1.525  & 1.546  & 1.525  & 6.176  & 3.983  & 1.524  & 1.546  & 1.524  & 6.176  & 3.983  & 6.979  & 8.002  & 4.670  \\
1     & 3    & 1.561  & 1.608  & 1.549  & 6.598  & 2.915  & 1.562  & 1.611  & 1.546  & 6.598  & 2.916  & 7.203  & 8.001  & 14.874  \\

1     & 4   & 2.060  & 2.106  & 2.057  & 6.602  & 3.226  & 2.072  & 2.107  & 2.053  & 6.602  & 3.226  & 7.203  & 8.001  & 14.475  \\

2     & 1     & 1.219  & 1.935  & 1.218  & 2.232  & 1.757  & 1.222  & 1.935  & 1.218  & 2.232  & 1.758  & 2.197  & 2.576  & 3.677  \\
2     & 2     & 1.303  & 1.893  & 1.303  & 2.232  & 1.794  & 1.306  & 1.893  & 1.305  & 2.232  & 1.794  & 2.197  & 2.576  & 3.631  \\
2     & 3   & 1.285  & 2.592  & 1.304  & 2.261  & 1.690  & 1.285  & 2.592  & 1.304  & 2.261  & 1.690  & 2.434  & 2.576  & 3.930  \\

2     & 4   & 1.390  & 2.473  & 1.403  & 2.277  & 1.747  & 1.390  & 2.473  & 1.404  & 2.277  & 1.746  & 2.279  & 2.576  & 3.883  \\

3     & 1     & 21.916  & 24.553  & 20.617  & 25.237  & 18.567  & 21.975  & 24.552  & 20.619  & 25.237  & 18.573  & 29.523  & 28.408  & 63.842  \\
3     & 2     & 21.034  & 24.100  & 20.617  & 25.237  & 21.182  & 21.039  & 24.100  & 20.619  & 25.237  & 21.187  & 29.523  & 28.816  & 62.478  \\
3     & 3   & 483.129  & 397.007  & 22.599  & 28.660  & 25.719  & 484.323  & 397.040  & 22.603  & 28.660  & 25.725  & 31.804  & 31.514  & 85.041  \\

3     & 4     & 217.062  & 64.592  & 22.946  & 28.660  & 25.901  & 218.677  & 64.590  & 22.950  & 28.660  & 25.906  & 31.804  & 31.514  & 82.431  \\

4     & 1     & 5.651  & 6.627  & 5.494  & 6.169  & 5.742  & 5.653  & 6.627  & 5.496  & 6.169  & 5.744  & 6.392  & 6.358  & 13.230  \\
4     & 2     & 5.621  & 6.095  & 5.547  & 6.169  & 5.625  & 5.621  & 6.095  & 5.548  & 6.169  & 5.626  & 6.392  & 6.910  & 12.489  \\

4     & 3    & 60.014  & 153.772  & 5.940  & 6.759  & 6.028  & 59.144  & 153.772  & 5.941  & 6.759  & 6.030  & 6.982  & 7.318  & 15.286  \\
4     & 4    & 57.961  & 26.775  & 6.106  & 6.759  & 6.243  & 58.566  & 26.775  & 6.106  & 6.759  & 6.245  & 6.982  & 7.318  & 14.712  \\
\bottomrule
\end{tabular}%
\end{sidewaystable}

\begin{sidewaystable}[H]
\caption{Simulation results: Percentage of zeros in the combination weights (\%).}
\label{tab:weight_precent_zero}%
\footnotesize
\centering
\begin{tabular}{cc|*{14}{c}}
\toprule
Case & Set & $\hat{\w}^\calA_{\reg}$ & $\hat{\w}^\calB_{\reg}$ & $\hat{\w}^\calC_{\reg}$ & $\hat{\w}^\calD_{\reg}$ & $\hat{\w}^\calE_{\reg}$ & $\hat{\w}^\calA_{\ma}$ & $\hat{\w}^\calB_{\ma}$ & $\hat{\w}^\calC_{\ma}$ & $\hat{\w}^\calD_{\ma}$ & $\hat{\w}^\calE_{\ma}$ & $\hat{\w}^\calA_{\pf}$ & $\hat{\w}^\calB_{\pf}$ & $\hat{\w}^\calE_{\eig}$ \\
\midrule
1     & 1     & 0.00 & 0.00 &   9.09 &  45.45 &  18.18  & 0.00  & 0.00  &   9.09  &  45.45  & 18.18  &  90.91  &  72.73  & 0.00  \\
1     & 2     & 0.00 & 0.00 &  10.00 &  40.00 &  20.00  & 0.00  & 0.00  &  20.00  &  40.00  & 20.00  &  90.00  &  70.00  & 0.00  \\
1     & 3     & 0.00 & 0.00 &  23.68 &  63.16 &  15.79  & 0.00  & 0.00  &  28.95  &  63.16  &  0.00  &  97.37  &  68.42  & 0.00  \\
1     & 4     & 5.56 & 0.00 &  16.67 &  61.11 &  22.22  & 2.78  & 0.00  &  25.00  &  61.11  & 22.22  &  97.22  &  66.67  & 0.00  \\
2     & 1     & 0.00 & 0.00 &  18.18 &  45.45 &   0.00  & 0.00  & 0.00  &  27.27  &  45.45  &  0.00  &  81.82  &  90.91  & 0.00  \\
2     & 2     & 0.00 & 0.00 &  20.00 &  40.00 &   0.00  & 0.00  & 0.00  &  30.00  &  40.00  &  0.00  &  80.00  &  90.00  & 0.00  \\
2     & 3     & 0.00 & 0.00 &  42.11 &  86.84 &   0.00  & 0.00  & 0.00  &  42.11  &  86.84  & 2.63   &  97.37  &  89.47  & 0.00  \\
2     & 4     & 0.00 & 0.00 &  33.33 &  88.89 &   0.00  & 2.78  & 0.00  &  38.89  &  88.89  & 0.00   &  94.44  &  88.89  & 0.00  \\
3     & 1     & 0.00 & 0.00 &  45.45 &  54.55 &   0.00  & 0.00  & 0.00  &  45.45  &  54.55  & 0.00   &  90.91  &  90.91  & 0.00  \\
3     & 2     & 0.00 & 0.00 &  40.00 &  50.00 &   0.00  & 0.00  & 0.00  &  40.00  &  50.00  & 0.00   &  90.00  &  90.00  & 0.00  \\
3     & 3     & 0.00 & 0.00 &  73.68 &  89.47 &   0.00  & 0.00  & 0.00  &  73.68  &  89.47  & 0.00   &  97.37  &  97.37  & 0.00  \\
3     & 4     & 0.00 & 2.78 &  75.00 &  88.89 &   0.00  & 0.00  & 2.78  &  75.00  &  88.89  & 0.00   &  97.22  &  97.22  & 0.00  \\
4     & 1     & 0.00 & 0.00 &  36.36 &  54.55 &   0.00  & 0.00  & 0.00  &  36.36  &  54.55  & 0.00   &  90.91  &  90.91  & 0.00  \\
4     & 2     & 0.00 & 0.00 &  40.00 &  50.00 &   0.00  & 0.00  & 0.00  &  40.00  &  50.00  & 0.00   &  90.00  &  90.00  & 0.00  \\
4     & 3     & 0.00 & 0.00 &  73.68 &  78.95 &   0.00  & 0.00  & 0.00  &  73.68  &  78.95  & 0.00   &  97.37  &  97.37  & 0.00  \\
4     & 4     & 0.00 & 0.00 &  75.00 &  77.78 &   0.00  & 0.00  & 0.00  &  75.00  &  77.78  & 0.00   &  97.22  &  97.22  & 0.00  \\
\bottomrule
\end{tabular}%
\end{sidewaystable}

\section{Conclusion}\label{sec:conclusion}

In this paper, we highlighted the importance of the weight constraints or the region used to perform the optimization to find the optimal weights in forecast or model averaging. The constraints affect the properties of the combination and deserve attention in theoretical and applied papers.
Our suggestion is to avoid the default selection based on the convention and shift toward a more conscious approach that focuses on desired characteristics.
Specifically, if the in-sample fit is the main target, then unconstrained weights with the same objective function as the target criterion (e.g., SSR) leads to the best fit, while more constraints are typically associated with worse in-sample fit.
As a tradeoff, if the out-of-sample MSFE is the objective, then imposing more regulations and constraints often helps to reduce the variance and narrow down the upper bound of the combination MSFE.
The sum-up-to-unity constraint is a requisite when the focus is to guarantee empirical unbiasedness, while the positivity constraint is particularly useful if researchers would like to combine forecasts with only a small number of candidates, which may facilitate interpretation and reduce uncertainty; see also \citet{Radchenko2023} for a more detailed discussion on the role and treatment of negative weights.
Our discussion is based on several widely used objective functions, but more research is needed for recently proposed weights, for example, \citet{qian2022}, \citet{Christopher2024}, \citet{shi:su:xie:2022}, etc.

\bibliography{refs}
\bibliographystyle{abbrvnat}

\section*{Appendix}
\subsection*{Appendix A}\label{sec:appendix_A}
\begin{proof}{Proposition~\ref{pro:var}}
By constructions and with some algebra, we have:
\begin{align}
\Var_*(\hat{y}^\calA_{\reg,T+1})
&=\Var_*(\f_{T+1}\tp\hat{\w}^\calA_\reg)\notag\\
&=\f_{T+1}\tp\Var_*(\hat{\w}^\calA_\reg)\f_{T+1}\notag\\
&=\f_{T+1}\tp\Var_*\{(\F\tp\F)^{-1}\F\tp\y\}\f_{T+1}\notag\\
&=\sigma^2\f_{T+1}\tp(\F\tp\F)^{-1}\f_{T+1}\label{eq:VarA}
\end{align}
and
\begin{align}
\Var_*(\hat{y}^{\calA'}_{\reg,T+1})
&=\Var_*(\hat{\delta}_0+\f_{T+1}\tp\hat{\w}^{\calA'}_\reg)\notag\\
&=\tilde{\f}_{T+1}\tp\Var_*\left\{(\hat{\delta}_0, \hat{\w}^{{\calA'}\tp}_\reg)\tp\right\}\tilde{\f}_{T+1}\notag\\
&=\tilde{\f}_{T+1}\tp\Var_*\{(\tilde{\F}\tp\tilde{\F})^{-1}\tilde{\F}\tp\y\}\tilde{\f}_{T+1}\notag\\
&=\sigma^2\tilde{\f}_{T+1}\tp(\tilde{\F}\tp\tilde{\F})^{-1}\tilde{\f}_{T+1}\notag\\
&=\sigma^2\tilde{\f}_{T+1}\tp
\begin{pmatrix}
\one\tp\one & \one\tp\F \\
\F\tp\one &  \F\tp\F
\end{pmatrix}^{-1}
\tilde{\f}_{T+1}\notag\\
&=\sigma^2\tilde{\f}_{T+1}\tp
\begin{pmatrix}
\theta^{-1} &  -\theta^{-1}\one\tp\F(\F\tp\F)^{-1} \\
-\theta^{-1}(\F\tp\F)^{-1}\F\tp\one & (\F\tp\F)^{-1}+\theta^{-1}(\F\tp\F)^{-1}\F\tp\one\one\tp\F(\F\tp\F)^{-1}
\end{pmatrix}
\tilde{\f}_{T+1}\notag\\
&=\sigma^2\{\theta^{-1}-2\theta^{-1}\one\tp\F(\F\tp\F)^{-1}\f_{T+1}+\f_{T+1}\tp(\F\tp\F)^{-1}\f_{T+1}\notag\\
&\quad+\theta^{-1}\f_{T+1}\tp(\F\tp\F)^{-1}\F\tp\one\one\tp\F(\F\tp\F)^{-1}\f_{T+1}\}\notag\\
&=\sigma^2\left\{\theta^{-1}-2\theta^{-1}\beta+\f_{T+1}\tp(\F\tp\F)^{-1}\f_{T+1}+\theta^{-1}\beta^2\right\}\notag\\
& =\sigma^2\theta^{-1}(1-\beta)^2+\sigma^2\f_{T+1}\tp(\F\tp\F)^{-1}\f_{T+1}\notag\\
&\geq \Var_*(\hat{y}^\calA_{\reg,T+1}),\label{eq:VarB}
\end{align}
where $\tilde{\f}_{T+1}=(1, \f_{T+1}\tp)\tp$, $\tilde{\F}=(\one, \F)$, $\theta=n-\one\tp\F(\F\tp\F)^{-1}\F\tp\one$ and $\beta=\f_{T+1}\tp(\F\tp\F)^{-1}\F\tp\one$.

For $\hat{\w}^\calB_\reg$, we have:
\begin{align}
&\quad\Var_*(\hat{y}_{\reg,T+1}^\calB)\notag\\
&=\Var_*\{\f_{T+1}\tp\hat{\w}^\calA_\reg-\hat{\rho}_0\f_{T+1}\tp(\F\tp\F)^{-1}\one\}\notag\\
&=\f_{T+1}\tp\Var_*\{(\F\tp\F)^{-1}\F\tp\y-\phi^{-1}(\one\tp\hat{\w}^\calA_\reg-1)(\F\tp\F)^{-1}\one\}\f_{T+1}\notag\\
&=\f_{T+1}\tp\Var_*[(\F\tp\F)^{-1}\{\F\tp\y-\phi^{-1}(\one\tp(\F\tp\F)^{-1}\F\tp\y
-1)\one\}]\f_{T+1}\notag\\
&=\f_{T+1}\tp\Var_*[(\F\tp\F)^{-1}\{\F\tp-\phi^{-1}\one\one\tp(\F\tp\F)^{-1}\F\tp\}\y]\f_{T+1}\notag\\
&=\sigma^2\f_{T+1}\tp(\F\tp\F)^{-1}\{\I-\phi^{-1}\one\one\tp(\F\tp\F)^{-1}\}\F\tp\F
\{\I-\phi^{-1}(\F\tp\F)^{-1}\one\one\tp\}(\F\tp\F)^{-1}\f_{T+1}\notag\\
&=\sigma^2\f_{T+1}\tp(\F\tp\F)^{-1}\{\F\tp\F-\phi^{-1}\one\one\tp\}\{\I-\phi^{-1}(\F\tp\F)^{-1}\one\one\tp\}(\F\tp\F)^{-1}\f_{T+1}\notag\\
&=\sigma^2\f_{T+1}\tp(\F\tp\F)^{-1}\{\F\tp\F-\phi^{-1}\J_n\}\{\I-\phi^{-1}(\F\tp\F)^{-1}\J_n\}(\F\tp\F)^{-1}\f_{T+1}\notag\\
&=\sigma^2\f_{T+1}\tp(\F\tp\F)^{-1}\{\F\tp\F-2\phi^{-1}\J_n+\phi^{-1}\J_n\}(\F\tp\F)^{-1}\f_{T+1}\notag\\
&=\sigma^2\f_{T+1}\tp(\F\tp\F)^{-1}\f_{T+1}-\phi^{-1}\sigma^2
\{\f_{T+1}\tp(\F\tp\F)^{-1}\one\}^2\notag\\
&\leq \sigma^2\f_{T+1}\tp(\F\tp\F)^{-1}\f_{T+1}=\Var_*(\hat{y}^\calA_{\reg,T+1})\label{eq:VarC}
\end{align}
where $\phi=\one\tp(\F\tp\F)^{-1}\one$, $\J_n=\one\otimes\one\tp$.
From \eqref{eq:VarA}, \eqref{eq:VarB} and \eqref{eq:VarC}, we know the variance of the combined forecast is decreased by imposing the constraint $\one\tp\w=1$, and $\Var_*(\hat{y}_{\reg,T+1}^\calB)\leq\Var_*(\hat{y}^\calA_{\reg,T+1})\leq \Var_*(\hat{y}_{\reg,T+1}^{\calA'})$. The other optimal weights of CV model averaging in $\bbW^\calA$ and $\bbW^\calB$ have the same result and here we give the conclusion without proofs:
\begin{align*}
\Var_*(\hat{y}_{\cv,T+1}^{\calA})\geq \Var_*(\hat{y}_{\cv,T+1}^\calB).
\end{align*}

Next, for $\hat{\w}^\calC_Z\in\bbW^\calC$, where the subscript $Z$ represents ``reg'', ``ma'' and ``cv'', we have:
\begin{align*}
\Var_*(\hat{y}^\calC_{Z,T+1})&=\Var_*(\f_{T+1}\tp \hat{\w}_Z^\calC)\notag\\
&=\f_{T+1}\tp\Var_*(\hat{\w}_Z^\calC)\f_{T+1}\notag\\
&\leq \f_{T+1}\tp\Exp_*\{\hat{\w}_Z^\calC(\hat{\w}_Z^\calC)\tp\}\f_{T+1}\notag\\
&\leq \f_{T+1}\tp\Exp_*\left[\lambda_{\max}\{\hat{\w}_Z^\calC(\hat{\w}_Z^\calC)\tp\}\right]\f_{T+1}\notag\\
&= \f_{T+1}\tp\Exp_*(\|\hat{\w}_Z^\calC\|^2)\f_{T+1}\notag\\
&\leq S \f_{T+1}\tp\f_{T+1}.
\end{align*}

For $\hat{\w}^\calD_Z\in\bbW^\calD$, where the subscript $Z$ represents ``reg'', ``ma'' and ``cv'' ``pf'', we have:
\begin{align*}
\Var_*(\hat{y}^\calD_{Z,T+1})
&=\Var_*(\f_{T+1}\tp \hat{\w}^\calD_Z)\notag\\
&=\f_{T+1}\tp\Var_*(\hat{\w}^\calD_Z)\f_{T+1}\notag\\
&\leq \f_{T+1}\tp\Exp_*\{\hat{\w}^\calD_Z(\hat{\w}^\calD_Z)\tp\}\f_{T+1}\notag\\
&\leq \f_{T+1}\tp\Exp_*\left[\lambda_{\max}\{\hat{\w}^\calD_Z(\hat{\w}^\calD_Z)\tp\}\right]\f_{T+1}\notag\\
&= \f_{T+1}\tp(\Exp_*\|\hat{\w}^\calD_Z\|^2)\f_{T+1}\notag\\
&= \f_{T+1}\tp\f_{T+1}\Exp_*\left\{\sum_{s=1}^{S}(\hat{w}_{Z,s}^{\calD })^2\right\}\notag\\
&\leq  \f_{T+1}\tp\f_{T+1}\Exp_*\left\{\sum_{s=1}^{S}\hat{w}^\calD_{Z,s}\right\}^2\notag\\
&=\f_{T+1}\tp\f_{T+1},
\end{align*}
where the last inequality is because of the condition that $\hat{w}_s\geq 0$ for $s=1, 2, \cdots, S$ and the last equality is because $\hat{\w}\tp\one=1$.

Finally, for $\hat{\w}_\eig^\calE\in\bbW^\calE$, we have:
\begin{align*}
\Var_*(\hat{y}_{\eig,T+1}^\calE)
&=\Var_*(\f_{T+1}\tp \hat{\w}_\eig^\calE)\notag\\
&=\f_{T+1}\tp\Var_*(\hat{\w}_\eig^\calE)\f_{T+1}\notag\\
&\leq \f_{T+1}\tp\Exp_*\{\hat{\w}_\eig^\calE(\hat{\w}_\eig^\calE)\tp\}\f_{T+1}\notag\\
&\leq \f_{T+1}\tp\Exp_*\left[\lambda_{\max}\{\hat{\w}_\eig^\calE(\hat{\w}_\eig^\calE)\tp\}\right]\f_{T+1}\notag\\
&= \f_{T+1}\tp\Exp_*(\|\hat{\w}_\eig^\calE\|^2)\f_{T+1}\notag\\
&= \f_{T+1}\tp\f_{T+1}.
\end{align*}
\end{proof}

\subsection*{Appendix B}\label{sec:appendix_B}
\newcommand{\bLambda}{\mathbf{\Lambda}}
\begin{proof}{the uniqueness of $\hat{\w}^\calE_\reg$}
The weight $\hat{\w}^\calE_\reg$ is the optimal solution of the following optimization problem:
\begin{align}\label{eq:optimization}
\min_\w \|\y-\F\w\|^2\quad \text{s.t.}\quad \|\w\|^2=1.
\end{align}
We can construct the following objective function with a Lagrangian multiplier $\nu$ as
\begin{align*}
l(\w,\nu)=\|\y-\F\w\|^2-\nu(\w\tp\w-1).
\end{align*}
Set the first derivative of $l(\w,\nu)$ to be zero, then
\begin{align}
\w=(\F\tp\F-\nu)^{-1}\F\tp\y,\label{eq:deriv1}\\
\|(\F\tp\F-\nu)^{-1}\F\tp\y\|^2=1.\label{eq:deriv2}
\end{align}
Considering that $\F\tp\F$ is a positive definite matrix, according to the properties of symmetric matrices, there exists an orthogonal matrix $\Q$ such that $\F\tp\F=\Q\bLambda\Q^{-1}$, where $\bLambda=\diag(\lambda_1, \lambda_2, \ldots, \lambda_n)$  is a diagonal matrix consisting of eigenvectors of $\F\tp\F$, with $\lambda_i>0 \text{ for } i=1, \ldots, n$. Then,
from \eqref{eq:deriv2}, we have:
\begin{align*}
\|(\F\tp\F-\nu)^{-1}\F\tp\y\|^2&=\|(\Q\bLambda\Q^{-1}-\nu)^{-1}\F\tp\y\|^2\notag\\
&=\|\Q(\bLambda-\nu)^{-1}\Q^{-1}\F\tp\y\|^2\notag\\
&=\|\Q(\bLambda-\nu)^{-1}\Q^\top\F\tp\y\|^2\notag\\
&=\y\tp \F\Q (\bLambda-\nu)^{-1}\Q^\top \Q(\bLambda-\nu)^{-1}\Q^\top\F\tp\y\notag\\
&=\y\tp \F\Q (\bLambda-\nu)^{-2}\Q^\top\F\tp\y\notag\\
&=\tilde{\y}^\top (\bLambda-\nu)^{-2}\tilde{\y}\notag\\
&=\sum_{i=1}^{n}\frac{\tilde{y}_i^2}{(\lambda_i-\nu)^2},
\end{align*}
where $\tilde{\y}=\Q^\top\F\tp\y$.
The function $f(\nu)=\sum_{i=1}^{n}\frac{\tilde{y}_i^2}{(\lambda_i-\nu)^2}$ is monotonically decreasing in $(-\infty, \lambda_{\min})$ and increasing in $(\lambda_{\max},\infty)$, and $0=f(-\infty)<1<f(\lambda_{\min})=\infty, \infty=f(\lambda_{\max})>1>f(-\infty)=0$. Therefore, there exist two solutions $\nu_1\in (-\infty, \lambda_{\min})$ and $\nu_2\in (\lambda_{\max}, \infty)$.

Besides, if $\|(\F\tp\F-\nu)^{-1}\F\tp\y\|>1$ for $\nu\in (\lambda_{\min}, \lambda_{\max})$, which is equivalent to
\begin{align*}
\sum_{i=1}^{n}\frac{\tilde{y}_i^2}{(\lambda_i-\nu)^2}>1 \text{ for }\nu\in (\lambda_{\min}, \lambda_{\max}),
\end{align*}
then $\nu_1$ and $\nu_2$ are the only two solutions for \eqref{eq:deriv2}.
Furthermore, the Hessian matrix with respect to $\w$ can be obtained by
\begin{align*}
\frac{\partial^2 l(\w,\nu)}{\partial \w\partial \w\tp}=2(\F\tp\F-\nu\I)>0.
\end{align*}
It is positive definite when $\nu_1 \in(-\infty, \lambda_{\min})$, and negative definite when $\nu_2 \in(\lambda_{\max}, \infty).$ Hence, according to the convex optimization theories \citep{Stephen2004Convex}, $\w=(\F\tp\F-\nu_1)^{-1}\F\tp\y$ generated by \eqref{eq:deriv1} is the optimal solution for \eqref{eq:optimization} and it is unique.

\end{proof}

\end{document}